\documentclass[11pt]{article}
\usepackage[T1]{fontenc}
\usepackage{lmodern}
\usepackage[a4paper,margin=27mm]{geometry}
\usepackage{amsmath,amssymb,amsthm,mathtools,mathrsfs}
\usepackage{microtype,booktabs,array,enumitem}
\usepackage{xcolor,xurl}
\usepackage[colorlinks=true,linkcolor=blue!45!black,citecolor=blue!45!black,urlcolor=blue!45!black]{hyperref}
\usepackage{authblk}
\usepackage[nameinlink,noabbrev]{cleveref}
\numberwithin{equation}{section}
\allowdisplaybreaks[2]
\newcommand{\doi}[1]{\href{https://doi.org/#1}{\nolinkurl{doi:#1}}}
\newcommand{\R}{\mathbb R}
\newcommand{\C}{\mathbb C}

\newcommand{\PP}{\mathbb P}
\newcommand{\EE}{\mathbb E}
\newcommand{\dd}{\,\mathrm d}
\newcommand{\ii}{\mathrm i}

\newcommand{\pair}[2]{\left\langle #1,#2\right\rangle}
\newcommand{\dist}{\operatorname{dist}}
\newcommand{\supp}{\operatorname{supp}}
\newcommand{\tr}{\gamma}
\newcommand{\one}{\mathbf 1}
\newcommand{\cE}{\mathcal E}
\newcommand{\cH}{\mathcal H}
\newcommand{\cR}{\mathcal R}

\newcommand{\cL}{\mathcal L}

\newcommand{\cD}{\mathcal D}

\newcommand{\cX}{\mathcal X}
\newcommand{\cY}{\mathcal Y}

\newcommand{\fl}{\mathrm{FL}}
\newcommand{\eff}{\mathrm{eff}}
\newcommand{\refe}{\mathrm{ref}}
\newcommand{\av}{\mathrm{av}}
\newcommand{\cG}{\mathcal G}
\newtheorem{theorem}{Theorem}[section]
\newtheorem{proposition}[theorem]{Proposition}
\newtheorem{lemma}[theorem]{Lemma}
\newtheorem{corollary}[theorem]{Corollary}
\theoremstyle{definition}
\newtheorem{assumption}[theorem]{Assumption}

\theoremstyle{remark}
\newtheorem{remark}[theorem]{Remark}
\crefname{assumption}{Assumption}{Assumptions}
\crefname{remark}{Remark}{Remarks}
\title{The Calder\'on Problem with Randomly\\ Distributed Resonant Droplets}

\author{Anjana E. J. and Mourad Sini\thanks{The work of Anjana E. J. is supported by the DOC Fellowship, Austrian Academy of Sciences (\"{O}AW). The work of M. Sini is partially supported by the Austrian Science Fund (FWF): P 32660 and P 36942.}}

\affil{\small Johann Radon Institute (RICAM), Austrian Academy of Sciences, Linz, Austria (\texttt{anjana.jeevanand@ricam.oeaw.ac.at} and \texttt{mourad.sini@oeaw.ac.at}).}

\date{}
\begin{document}
\maketitle
\begin{abstract}

We study the recovery of an acoustic coefficient from boundary measurements
made after the injection of a collection of small, highly contrasting droplets. This enters into the framework of quantitative ultrasound imaging using contrast agents. The centers are randomly distributed according to a prescribed strictly positive density $\rho\in C^1(\overline\Omega)$. We choose real incident frequencies near an active Newtonian resonance and prove a point-interaction approximation, followed by a comparison with the
effective model $(-\Delta+P^2\rho-\omega^2n^2)u=0$. The constants are uniform as the amplitude $P^2$ increases, provided it is coupled to the size and distribution discrepancy of the droplets. A random rounding and displacement construction verifies the distribution assumption for each such given density. The effective Neumann-to-Dirichlet map is then linearized with a remainder of order $P^{-4}$. We derive reconstruction identities using complex geometrical optics solutions of the known reference operator $-\Delta+P^2\rho$. We also estimate the amplification of the boundary-map error in these identities and give a boundary integral correction for reconstruction from the full (and not the linearized) measured map.
\end{abstract}
\setcounter{tocdepth}{1}
\tableofcontents
\clearpage
\section{Introduction and statement of the results}
\label{sec:model}

\subsection{The inverse problem and related works}
The Calder\'on problem, we are interested in, asks whether an interior coefficient can be recovered from the corresponding boundary measurement map. In its original formulation, the coefficient is an electrical conductivity and the measurements are encoded by the Dirichlet-to-Neumann map \cite{Calderon}. The advancement of the problem began with the construction of complex geometrical optics (CGO) solutions, by Sylvester and Uhlmann, in \cite{SU}, which established global uniqueness for $n \ge 3$, see also the original work by Ramm \cite{Ramm-book}. Later on, Nachman \cite{Nachman88} and Novikov \cite{Novikov88} developed methods to compute the interior potential from the boundary measurements.

\noindent The initial CGO-based results led to a research direction to determine the minimal regularity conditions on the conductivity under which uniqueness holds. These assumptions have been relaxed from smooth coefficients down to conductivities in various fractional Sobolev and Besov spaces (see, \cite{Brown1996, Paivarinta, JerisonKenig, Brown2003, NguyenSpirin2014, Haberman2014}). Finally, uniqueness results for Lipschitz conductivities were established in \cite{HabermanTataru}, and further relaxations of conductivities are established in \cite{CaroRogers}. The two-dimensional problem requires different arguments, with global uniqueness proved using $\bar{\partial}$ methods \cite{Nachman96} and extended to $L^\infty$ conductivities using quasiconformal mappings \cite{AstalaPaivarinta}. Bukhgeim \cite{Bukhgeim} introduced an alternative framework using stationary phase techniques. For an account
of these developments, see \cite{Uhlmann,Salo}.

\noindent The stability estimates in \cite{Alessandrini} and the instability result in \cite{Mandache} describe an important difficulty of recovery from remote measurements. One way to remedy to this instability, and obtain additional information, is to perturb the medium by injecting contrast agents whose size, shape and material parameters are known. The resulting inverse problem uses the change in the measured field. When the agents are small and highly contrasting, their resonances can produce a leading-order perturbation
which contains information about the surrounding medium. This enters into the framework of imaging using contrast agents. In the recent years, we proposed a modeling and analysis of such approaches where the main idea is the use of the resonant behavior of the contrast agents in its different manifestations. In this direction, this approach has been analyzed for acoustic imaging with nearly resonating bubbles \cite{DGS}, for photo-acoustic imaging with dielectric particles \cite{GSdielectric}, and for plasmonic particles in the full Maxwell system \cite{GSMaxwell,GSsimultaneous}. Optical inversion using plasmonic contrast agents is considered in \cite{CGSoptical}. In the time domain, injected droplets have been used for recovering a wave source \cite{SW} and, subsequently, both the wave speed and the source \cite{SSW}. We refer to \cite{GSSsurvey} for an introduction to the mathematical analysis of imaging modalities using small contrast agents, and for a review of the resonance mechanisms and
their use in inverse problems. These works concern different measurement settings and different resonances. We use, in the current work, the acoustic droplet model and the boundary-map formulation of \cite{JST}.

\noindent In \cite{JST}, all the droplets are injected at once contrary to \cite{DGS, GSdielectric,GSMaxwell,GSsimultaneous, SW, SSW} where the agents were injected in isolation (i.e. one after another). A real incident frequency is chosen near a Newtonian resonance. The boundary map of the perturbed medium is shown to approach that of an effective medium with coefficient $\omega^2n^2-P^2$. The parameter $P^2>0$ can be made large by adjusting the detuning. This gives a coercive reference problem and allows the effective boundary map to be linearized before applying CGO solutions. 

\noindent Our purpose here is to replace the periodic distribution in \cite{JST} by a random distribution with a prescribed, possibly nonconstant, density. We formulate the configurations on a probability space and work on events on which the geometry and the empirical measure satisfy uniform bounds. For each realization in such an event, the proof is deterministic: we derive the Foldy--Lax system and compare it with a density-weighted integral equation. We do not use stochastic two-scale convergence. The dependence on $P$ is traced in each estimate. The use of geometrically constrained random configurations is also present in \cite{BBM}, where the authors establish resonant homogenization for dielectric rods with random positions, radii and permittivities through stochastic two-scale analysis. Here, we use the point-interaction approximation and quantitative empirical-measure estimates to derive an effective equation for randomly distributed acoustic droplets. Our motivation is to apply this approximation to the Calder\'on problem, with boundary-map estimates that retain the dependence on the tunable resonant amplitude.

\noindent The effective perturbation is $-P^2\rho(x)$. Thus, for a nonconstant spatial density, the reference operator is $-\Delta+P^2\rho(x)$ rather than $-\Delta+P^2$. We prove the corresponding resolvent estimates and linearization, and construct CGO solutions for this known reference operator. The boundary norms of these solutions grow with their complex phase. We therefore retain these norms when applying the linearization to measured data which quantify  to which extent the reconstruction is valid. Appendix~\ref{sec:correction} treats reconstruction from the full boundary map by a boundary integral equation, following the constructive approach of \cite{Nachman88,Novikov88} while the perturbation
estimates are related to those used in \cite{KR}.

\noindent Our reconstruction strategy relies on tuning the incident frequency near a common droplet resonance to generate a negative effective perturbation of controllably large magnitude. Random variations in the material parameters would shift the individual resonances and could weaken this collective response or introduce contributions of opposite signs. Preserving the required sign and amplitude would therefore require additional restrictions relating the resonance spread to the chosen detuning, which becomes smaller as the desired amplitude increases. We consequently keep the droplet size, shape and material parameters fixed and introduce randomness only in their centers and occupancies. This choice also reflects a natural idealization for contrast-agent imaging, namely the injected agents are prepared with prescribed and calibrated characteristics, whereas their individual locations after injection may not be controlled. Our analysis thus isolates spatial uncertainty while preserving the resonant mechanism used for reconstruction.
\bigskip

\noindent We would like to finish this introduction with the following statement. In the current work, we start from the original model $-\Delta-\omega^2n^2$, where $n$ is the unknown coefficient and $\omega$ is the operating incident frequency. Then after injecting the droplets, all at once, with a given distribution of their locations, and passing to the homogenized model, we have the two choices according to whether the used incident frequency is chosen close to the resonance from the left-side or the right-side. Accordingly, we have two outcome:
\begin{enumerate}
\item Using frequencies from the left-side generates an effective model with the background of the form $-\Delta+P^2\rho$. In this case, in terms of inverse problems, we can linearize the effective measurement map. The drawback is that we get limited band-frequency for the stability estimates with Fourier modes at most in the band $[0, \omega]$. 

\item Using incident frequencies at the right-side of the resonant frequency generate and effective background of the form $-\Delta-P^2\rho$. In this case, a linearization of the effective boundary map, in terms of $P$, is not possible. But the inversion of the non-linear map provides stable Fourier bands of the form $[0, \sqrt{\omega^2+P^2}]$ (in the case $\rho=1$), which suggest, in principle, and increase of the quantitative stability as $P$ grows. Such results will be reported in a forthcoming work.
\end{enumerate}
\noindent As a concluding message, using such contrast agents we can breakdown the two well known obstacles of the Calder\'on-problem, namely non-linearity and instability, in the following sense.
\begin{enumerate}
\item using left-side resonant frequencies, we can linearize the boundary map but we are limited in stability.

\item using right-side resonant frequencies, we cannot linearize the boundary map but we can improve the stability by increasing the stable Fourier-band. 
\end{enumerate}
Therefore, either we breakdown the non-linearity or the instability but not both. To simultaneously linearize and stabilize the inverse problem, one tentative is to choose the resonant incident frequencies, with a small damping, to be able to replace the effective potential $\pm P^2$ by a convenient complex parameter having convenient (large) real and imaginary parts that could allow linearization (with large imaginary part) and stability with (large) real part. But this understanding needs to be confirmed.
\bigskip

\noindent Let us emphasize that we are not suggesting to take the incident frequency large to improve stability, as it is suggested by Isakov \cite{Isakov-1} and others \cite{Isakov-2}, rather, we use the resonant effect to increase the effective background potential (i.e. $P^2$).

\subsection{The acoustic model}
Let $\Omega\subset\R^3$ be a bounded domain with $C^{2,\alpha}$ boundary,
where $0<\alpha<1$. The physical mass density is a constant $r_0>0$.
We write
\[
 n^2=\frac{r_0}{k_b}\in W^{1,\infty}(\Omega),\qquad
 0<n_*^2\le n^2\le n^{*2},
\]
where $k_b$ is the bulk modulus of the background. The notation $\rho$
will be reserved for the spatial density of the droplet centers. Set
\[
 \cX=H^{-1/2}(\partial\Omega),\qquad
 \cY=H^{1/2}(\partial\Omega),\qquad q_\omega=\omega^2n^2.
\]
Unless indicated otherwise, boundary duality pairings are bilinear,
obtained by complex-linear extension of the real duality. Conjugates
will be written explicitly in energy estimates.

We assume that there is a compact interval $I\Subset(0,\infty)$, with
$\omega_0$ in its interior, such that the Neumann realization of
$A_{0,\omega}=-\Delta-q_\omega$ satisfies \footnote{Tuning the droplets properties, i.e. its bulk for instance, allows us to choose $\omega_0$, see (\ref{Non-degenrate-moment}), and hence $I$, so that (\ref{eq:background-gap}) is satisfied.}
\begin{equation}\label{eq:background-gap}
 \sup_{\omega\in I}
 \|A_{0,\omega}^{-1}\|_{L^2(\Omega)\to L^2(\Omega)}\le C_I.
\end{equation}
Here $\omega_0$ is the droplet resonance specified below. For $f\in\cX$,
let $S^f$ solve
\begin{equation}\label{eq:background}
 A_{0,\omega}S^f=0\quad\text{in }\Omega,\qquad
 \partial_\nu S^f=f\quad\text{on }\partial\Omega,
 \qquad \Lambda_{0,\omega}f=\tr S^f.
\end{equation}
The assumption \eqref{eq:background-gap} gives
$\|S^f\|_{H^1(\Omega)}\le C\|f\|_{\cX}$, uniformly for $\omega\in I$.
A justification, including the kernel estimates used below, is given
in Appendix~\ref{sec:auxiliary}.

Let $G_\omega$ be the real symmetric Neumann Green kernel defined by
\begin{equation}\label{eq:green-sign}
 (-\Delta_x-q_\omega(x))G_\omega(x,y)=\delta_y,
 \qquad \partial_{\nu_x}G_\omega(x,y)=0.
\end{equation}
We use the decomposition
\[
 G_\omega(x,y)=\Phi_0(x,y)+H_\omega(x,y),\qquad
 \Phi_0(x,y)=\frac1{4\pi|x-y|}.
\]
The estimates needed in the sequel are, see Lemma \ref{lem:background-bounds},
\begin{align}
 |G_\omega(x,y)|&\le C|x-y|^{-1},&
 |\nabla_xG_\omega(x,y)|+|\nabla_yG_\omega(x,y)|
   &\le C|x-y|^{-2},\label{eq:green-bounds}\\
 |H_\omega(x,y)|&\le Cb^{-1},&
 |\nabla_xH_\omega(x,y)|+|\nabla_yH_\omega(x,y)|
   &\le Cb^{-2}.\label{eq:H-bounds}
\end{align}
In \eqref{eq:H-bounds}, $0<b\le b_1$, and $x,y$ belong to a ball of
radius comparable to $b$, whose distance from $\partial\Omega$ is at
least a fixed multiple of $b$. The value estimate holds on the diagonal;
the derivative estimate is understood off the diagonal, or as a
Lipschitz estimate. Constants are uniform for $\omega\in I$.

With the background framework and its associated resolvent estimates established, we now introduce the physical perturbation. The continuous background $\Omega$ is subsequently perturbed by injecting a large ensemble of $M$ randomly distributed, subwavelength acoustic droplets. The  geometric properties and resonant scaling of these individual inclusions are formalised as follows.

\begin{assumption}[The droplets and the resonance]\label{ass:micro}
Each droplet is of the form $D_j=z_j+aB$, where $B\subset\R^3$ is a
fixed bounded $C^{1,1}$ domain containing the origin. Its mass density
is $r_0$ and its bulk modulus is $k_0a^2$, with $k_0>0$. Let
\[
 (N_Bv)(\xi)=\int_B\frac{v(\eta)}{4\pi|\xi-\eta|}\dd\eta.
\]
Choose an isolated eigenvalue $\lambda_*>0$ of $N_B$ and denote its
orthogonal spectral projection by $\Pi_*$. We assume
\begin{equation}\label{Non-degenrate-moment}
 b_*:=\|\Pi_*\one\|_{L^2(B)}^2>0,
 \qquad \omega_0^2=\frac{k_0}{r_0\lambda_*}.
\end{equation}
Only the centers, and possibly their total number, are random.
\end{assumption}

Fix $1/3<h<1$ and introduce
\begin{equation}\label{eq:scales}
 w=a^{1-h},\qquad d=w^{1/3},\qquad t=P^2\ge t_0\ge1.
\end{equation}
The scale $d$ is deterministic. It provides a lower bound, up to a
constant, for the actual distance between distinct centers. We describe
the configurations on a probability space $(\Xi,\mathscr F,\PP)$ and set
$\varpi\in\Xi$ to denote a realization. For each $a$, $M(a,\varpi)$ denotes the total number of
droplets in realization $\varpi$. The map $M(a,\cdot)$ is
integer-valued and measurable, and $z_j$ is measurable on
$\{M(a,\cdot)\ge j\}$.

\begin{assumption}[Separation and boundary clearance]\label{ass:hardcore}
There are measurable events $\cE_a\subset\Xi$ and deterministic
constants $c_0,c_b,C_M>0$ such that, for $\varpi\in\cE_a$,
\[
 \min_{i\ne j}|z_i-z_j|\ge c_0d,\qquad
 \min_j\dist(z_j,\partial\Omega)\ge c_bd,\qquad
 wM\le C_M.
\]
The events and constants do not depend on $P$ or on the boundary input.
\end{assumption}

Since $a/d\to0$, this assumption implies, for sufficiently small $a$,
that the droplets are disjoint and contained in $\Omega$. Their union
will be denoted by $D$.

\begin{assumption}[The distribution of the centers]\label{ass:LLN}
Let $\rho\in C^1(\overline\Omega)$ be known and satisfy
$0<\rho_*\le\rho\le\rho^*$. Define
\[
 \mu_a(\varpi):=w\sum_{j=1}^{M(a,\varpi)}\delta_{z_j(\varpi)},\qquad
 \mathfrak D_a(\varpi)=:
 \sup_{\|\varphi\|_{C^1(\overline\Omega)}\le1}
 \left|\int\varphi\dd\mu_a-\int_\Omega\rho\varphi\dd x\right|,
\]
where $\|\varphi\|_{C^1}=\|\varphi\|_\infty+\|\nabla\varphi\|_\infty$.
There are deterministic quantities $D_a\to0$ and $p_a\to0$ such that
\begin{equation}\label{eq:common-event}
 \PP(\cE_a)\ge1-p_a,\qquad
 \mathfrak D_a(\varpi)\le D_a\quad\text{for }\varpi\in\cE_a.
\end{equation}
\end{assumption}

The supremum in this assumption can be taken over a countable dense
subset of the $C^1$ unit ball, so that $\mathfrak D_a$ is measurable.
This is a common-event estimate, not a law of large numbers for one
fixed test function. We shall apply it to regularized Green kernels and
to products of functions depending on $P$. It also determines the
normalization of the total number of droplets:
\begin{equation}\label{eq:mass-normalization}
 \left|wM-\int_\Omega\rho\dd x\right|\le\mathfrak D_a.
\end{equation}
In particular, no normalization $\int_\Omega\rho=1$ is imposed.
Section~\ref{sec:random} constructs random configurations satisfying
Assumptions~\ref{ass:hardcore} and \ref{ass:LLN} for every fixed density
in this class.

\subsection{The near-resonant frequency and the effective equation}
We choose the real frequency
\begin{equation}\label{eq:tuning}
 \omega_{a,P}^2=\omega_0^2
 \left(1+\frac{b_*}{\lambda_*P^2}a^h\right),\qquad
 c_*(P)=-\frac{k_0b_*}{\lambda_*P^2}<0.
\end{equation}
For fixed $P$, this is the detuning in \cite[equation (1.6)]{JST}.
The choice of its sign will give a negative scattering coefficient.
For $a$ sufficiently small, the frequencies in \eqref{eq:tuning} belong
to $I$, uniformly for $P\ge P_0$. In every comparison below, all the
boundary maps are evaluated at this same frequency. We often write
$\omega=\omega_{a,P}$ and suppress it in $G_\omega$.

The exact acoustic field is the solution of
\begin{equation}\label{eq:exact-PDE}
 \begin{cases}
 \displaystyle
 \left[-\Delta-q_\omega(1-\one_D)
       -\frac{\omega^2r_0}{k_0a^2}\one_D\right]u_a^f=0
       &\text{in }\Omega,\\
 \partial_\nu u_a^f=f&\text{on }\partial\Omega.
 \end{cases}
\end{equation}
We set $\Lambda_{a,P}(\varpi)f=\tr u_a^f$. The constant mass density
makes \eqref{eq:exact-PDE} a scalar equation with a zeroth-order
contrast (i.e. no contrast occurs in its principal part). The effective field
and its boundary map are defined by
\begin{equation}\label{eq:effective-PDE}
 \begin{cases}
 (-\Delta+P^2\rho-q_\omega)u_{\rho,P}^f=0&\text{in }\Omega,\\
 \partial_\nu u_{\rho,P}^f=f&\text{on }\partial\Omega,
 \end{cases}
 \qquad \Lambda_{\rho,P}^{\eff}f=\tr u_{\rho,P}^f.
\end{equation}

\begin{theorem}\label{thm:main}
Suppose that Assumptions~\ref{ass:micro}--\ref{ass:LLN} hold. Put
\begin{equation}\label{eq:epsilon}
 \epsilon_a=(D_a+d)^{1/5},\qquad
 E(a,P)=P^2a^h+P^4d+P\epsilon_a.
\end{equation}
There are constants $a_0,c_{\mathrm{adm}},C>0$ and $P_0\ge1$ such that,
if
\begin{equation}\label{eq:joint-regime}
 0<a<a_0,\qquad P\ge P_0,\qquad P\epsilon_a\le c_{\mathrm{adm}},
\end{equation}
then \eqref{eq:exact-PDE} is uniquely solvable for every
$\varpi\in\cE_a$, and
\begin{equation}\label{eq:main-error}
 \|\Lambda_{a,P}(\varpi)-\Lambda_{\rho,P}^{\eff}\|_{\cL(\cX,\cY)}
 \le C\bigl(P^2a^h+P^4d+P\epsilon_a\bigr).
\end{equation}
Moreover, the Foldy--Lax matrix
\begin{equation}\label{eq:ideal-FL}
 (B_{a,P}Y)_i=Y_i+P^2w\sum_{j\ne i}G_\omega(z_i,z_j)Y_j
\end{equation}
is invertible and satisfies
\begin{equation}\label{eq:uniform-FL-inverse}
 \|B_{a,P}^{-1}\|_{\ell_w^2\to\ell_w^2}\le C,
 \qquad \|Y\|_{\ell_w^2}^2=w\sum_i|Y_i|^2.
\end{equation}
The constants are independent of $a,P,M,\varpi$ and of the boundary
inputs. The assertions hold on the same event for all the parameters
in \eqref{eq:joint-regime}.
\end{theorem}

\noindent Here and below, the constants may depend on $\Omega,B,h,k_0,r_0$, the
spectral gap of the chosen eigenvalue, the bounds on $n^2$, the
background resolvent bound, the geometric constants in
Assumption~\ref{ass:hardcore}, and
$\rho_*^{-1},\rho^*,\|\rho\|_{C^1}$. Uniformity over a class of densities
means that these bounds are common to the class. The strict positivity
of $\rho$ is used in the large-$P$ coercivity and linearization estimates.

The proof is given in Sections~\ref{sec:resonance}--\ref{sec:boundary}.
We first compute the scattering coefficient, then establish the
continuous and discrete inverse estimates, and finally derive the
point-interaction system and compare the boundary maps. The weighted
norm in \eqref{eq:uniform-FL-inverse} makes these comparisons independent
of the growing number of droplets.

Let us emphasize the difference between the present argument and the one
in \cite{JST}. In \cite[Assumption~1]{JST}, the centers are deterministic:
there is one droplet in each translated cube of volume $a^{1-h}=w$, the
array is periodic in the interior of the domain, and the boundary layer is
handled deterministically. In that setting the discrete measure
\[
 \mu_a^{\rm Periodic}=w\sum_j\delta_{z_j}
\]
is a deterministic Riemann-sum measure; after the normalization used in
\cite{JST}, its limit is Lebesgue measure. Consequently, the passage from
the Foldy--Lax system to the continuous equation is based on cellwise
quadrature for a fixed configuration.

Here neither the number of centers in a cell nor their positions are
fixed. The basic object is the random empirical measure
\[
 \mu_a(\varpi)=w\sum_{j=1}^{M(a,\varpi)}\delta_{z_j(\varpi)},
\]
and the deterministic periodic quadrature condition of \cite{JST} is
replaced by the common-event estimate
\begin{equation}\label{eq:JST-random-comparison}
 \PP\left\{
   \sup_{\|\varphi\|_{C^1(\overline\Omega)}\le1}
   \left|\int\varphi\,\dd\mu_a-
          \int_\Omega\rho\varphi\,\dd x\right|
   \le D_a
 \right\}\ge1-p_a,
\end{equation}
combined with the hard-core separation and boundary-clearance conditions in Assumption \ref{ass:hardcore}.
Thus the randomness is isolated in the construction of an event
$\cE_a$ of large probability. Once $\varpi\in\cE_a$ is fixed, every
estimate in Sections~\ref{sec:distribution}--\ref{sec:boundary} is a
deterministic estimate, with constants independent of the realization.
In particular, no expectation of the field, no stochastic corrector and
no stochastic two-scale limit enters the Foldy--Lax analysis. This is also
different in purpose from the stochastic homogenization framework of
\cite{BBM}, namely we use a probability space to quantify admissible random
configurations, but the effective-medium comparison is proved
realization by realization on $\cE_a$.

There are three further consequences of this change. First, the event in
\eqref{eq:JST-random-comparison} controls the whole $C^1$ unit ball on one
and the same realization. It is stronger than convergence for a single
fixed test function and can therefore be applied to the regularized Green
kernels and to $P$-dependent products used later. Second, the limiting
coefficient is $P^2\rho(x)$, rather than the spatially constant effective
coefficient produced by the periodic unit-density array in \cite{JST}.
Third, since $\cE_a$ is independent of $P$ and of the boundary inputs, the
deterministic estimates obtained on $\cE_a$ hold simultaneously for all
admissible $P$ and all boundary data. This produces both the high-
probability statement below and, along summable sequences of scales, an
almost-sure statement by Borel--Cantelli.

\begin{corollary}\label{cor:probability}
For a deterministic choice $P=P(a)$ satisfying \eqref{eq:joint-regime},
\[
 \PP\left(
 \|\Lambda_{a,P}-\Lambda_{\rho,P}^{\eff}\|\le CE(a,P)
 \right)\ge1-p_a.
\]
If $E(a,P(a))\to0$, the difference of these maps converges to zero in
probability. If $a_k\downarrow0$ and $\sum_kp_{a_k}<\infty$, there is
a measurable set $\Xi_{\mathrm{typ}}\subset\Xi$ of probability one
such that \eqref{eq:main-error} holds for every
$\varpi\in\Xi_{\mathrm{typ}}$ and all sufficiently large $k$,
simultaneously for all admissible $P$.
\end{corollary}
\begin{proof}
The boundary-map error is measurable. To see this, fix $a,P,M$ and
regard the centers as parameters. On the fixed space $L^2(\Omega)$,
the Lippmann--Schwinger operator is
\[
 T_D=A_{0,\omega}^{-1}
       \bigl[(\omega^2r_0/(k_0a^2)-q_\omega)\one_D\,\cdot\bigr].
\]
Translations of the droplets are continuous in measure, and
$\sup_y\|G_\omega(\cdot,y)\|_2\le C$. Hence $T_D$ is continuous in
Hilbert--Schmidt norm. On the open set where $I-T_D$ is invertible,
the field is continuous in $L^2$ for each fixed boundary input. The
volume representation and $A_{0,\omega}^{-1}:L^2\to H^2$ then give
continuity of its boundary trace. Countable dense subsets of the unit
balls of $\cX$ and $\cY^*$ determine the operator norm, proving its
measurability. Taking the countable union over $M$ gives the same
conclusion for the random configuration.

Put
\[
 \mathcal A_a(P)=\left\{\varpi:\
 \|\Lambda_{a,P}(\varpi)-\Lambda_{\rho,P}^{\eff}\|
 \le CE(a,P)\right\}.
\]
Theorem~\ref{thm:main} states precisely that
$\cE_a\subset\mathcal A_a(P)$ for every admissible $P$, with the same
constant $C$. Therefore
\[
 \PP(\mathcal A_a(P))\ge\PP(\cE_a)\ge1-p_a,
\]
which is the first assertion.

Now let $P=P(a)$ be deterministic and assume $E(a,P(a))\to0$. Fix
$\eta>0$. There exists $a_\eta>0$ such that
$CE(a,P(a))\le\eta$ whenever $a<a_\eta$. For these values of $a$,
\[
 \left\{\|\Lambda_{a,P(a)}-\Lambda_{\rho,P(a)}^{\eff}\|>\eta\right\}
 \subset \cE_a^c.
\]
Hence
\[
 \PP\left\{\|\Lambda_{a,P(a)}-\Lambda_{\rho,P(a)}^{\eff}\|>\eta\right\}
 \le p_a\longrightarrow0,
\]
which is exactly convergence in probability.

For the almost-sure statement, set $F_k=\cE_{a_k}^c$. By
\eqref{eq:common-event},
\[
 \sum_{k=1}^\infty\PP(F_k)
 \le\sum_{k=1}^\infty p_{a_k}<\infty.
\]
The first Borel--Cantelli lemma, which does not require independence,
gives $\PP(\limsup_kF_k)=0$. Equivalently,
\begin{equation}\label{eq:typical-set-corollary}
 \Xi_{\mathrm{typ}}
 :=(\limsup_kF_k)^c
 =\bigcup_{K=1}^{\infty}\bigcap_{k\ge K}\cE_{a_k},
 \qquad \PP(\Xi_{\mathrm{typ}})=1.
\end{equation}
Thus, for every $\varpi\in\Xi_{\mathrm{typ}}$, there exists a finite
$K(\varpi)$ such that $\varpi\in\cE_{a_k}$ for all
$k\ge K(\varpi)$. Applying Theorem~\ref{thm:main} at each of those
scales gives \eqref{eq:main-error}. Since membership in $\cE_{a_k}$
does not involve $P$, and the theorem is uniform in $P$, the same
$K(\varpi)$ works simultaneously for every $P$ satisfying
\eqref{eq:joint-regime} at the corresponding scale, and also for every
boundary input.
\end{proof}
\section{The scattering coefficient}
\label{sec:resonance}

We begin with the single-droplet computation. It is the spectral
calculation of \cite[Appendix A.5]{JST}, with the detuning parameter
written in terms of $P$. Set
\begin{equation}\label{eq:theta}
 \theta=\frac{\omega_{a,P}^2r_0}{k_0}
 =\frac1{\lambda_*}\left(1+\frac{b_*a^h}{\lambda_*t}\right),
 \qquad \kappa_a=\theta a^{-2},\qquad t=P^2.
\end{equation}
On $L^2(B)$, define
\begin{equation}\label{eq:local-defs}
 R_a=(I-\theta N_B)^{-1},\qquad
 \psi_a=R_a\one,\qquad m_a=\int_B\psi_a,\qquad
 \alpha_a=a\theta m_a.
\end{equation}
The operators and quantities in \eqref{eq:local-defs} also depend on
$P$; this dependence is suppressed in the notation.

\begin{lemma}\label{lem:local}
For $a$ sufficiently small and $t\ge t_0\ge1$, the inverse $R_a$ exists.
There are constants independent of $a$ and $t$ such that
\begin{align}
 m_a&=-\lambda_*t a^{-h}+O(1),&
 \alpha_a&=-t a^{1-h}+a r_a,\quad |r_a|\le C,
 \label{eq:alpha-uniform}\\
 \|R_a\|_{L^2(B)\to L^2(B)}+\|\psi_a\|_{L^2(B)}
   &\le Ct a^{-h},&
 c t a^{-h}&\le |m_a|\le Ct a^{-h}.
 \label{eq:local-resolvent}
\end{align}
In particular, $\alpha_a<0$. The functional
\begin{equation}\label{eq:ell-functional}
 \ell_a(v)=m_a^{-1}\int_B R_av
          =\int_B\frac{\psi_a}{m_a}v
\end{equation}
is bounded on $L^2(B)$ with norm at most $C$, and $\ell_a(\one)=1$.
\end{lemma}
\begin{proof}
The Newtonian operator is compact and self-adjoint on $L^2(B)$.
Since $\lambda_*$ is isolated, there is $\delta_*>0$ such that
\[
 \left|1-\frac{\lambda}{\lambda_*}\right|\ge\delta_*
 \quad\text{for }\lambda\in\sigma(N_B)\setminus\{\lambda_*\}.
\]
The possible accumulation point $0$ causes no difficulty in this
inequality. By \eqref{eq:theta}, $\theta\to\lambda_*^{-1}$ uniformly
for $t\ge t_0$. After reducing $a_0$, we therefore have
$|1-\theta\lambda|\ge\delta_*/2$ on the complementary spectrum.
The spectral decomposition gives
\begin{equation}\label{eq:resolvent-splitting}
 R_a=\frac{\Pi_*}{1-\theta\lambda_*}+R_a^\perp,
 \qquad R_a^\perp=(I-\theta N_B)^{-1}(I-\Pi_*),
 \qquad \|R_a^\perp\|\le\frac2{\delta_*}.
\end{equation}
From (\ref{eq:theta}), we have,
\[
 1-\theta\lambda_*=-\frac{b_*a^h}{\lambda_*t}
\]
thus the denominator is nonzero. Taking the inner product of
\eqref{eq:resolvent-splitting} with $\one$ yields
\begin{equation}\label{eq:mass-spectral}
 m_a=\frac{b_*}{1-\theta\lambda_*}
       +\pair{R_a^\perp\one}{\one}_{L^2(B)}
     =-\lambda_*t a^{-h}+\beta_a,
 \qquad |\beta_a|\le C.
\end{equation}
Multiplication by $a\theta$ now gives
\begin{align*}
 \alpha_a
 &=a\frac1{\lambda_*}
       \left(1+\frac{b_*a^h}{\lambda_*t}\right)
       \left(-\lambda_*t a^{-h}+\beta_a\right)\\
 &=-t a^{1-h}
       +a\left(-\frac{b_*}{\lambda_*}+\theta\beta_a\right).
\end{align*}
The term in parentheses is uniformly bounded, because $\theta$ is
uniformly bounded. This proves \eqref{eq:alpha-uniform} without placing
any $t$-dependence in its remainder.

The operator estimate in \eqref{eq:local-resolvent} follows directly
from \eqref{eq:resolvent-splitting}. Since $t a^{-h}\to\infty$
uniformly for $t\ge t_0$, the bounded term $\beta_a$ can be absorbed
in \eqref{eq:mass-spectral}, proving the two-sided estimate for $m_a$.
The same observation shows that the main term in $\alpha_a$ dominates
$a r_a$, so that $\alpha_a<0$.

Finally, $R_a$ is real and self-adjoint. Consequently
$\int_B R_av=\int_B\psi_a v$, first for real $v$ and then by complex
linearity. By Cauchy--Schwarz and \eqref{eq:local-resolvent},
\[
 |\ell_a(v)|\le
 \frac{\|\psi_a\|_{L^2(B)}}{|m_a|}\|v\|_{L^2(B)}
 \le C\|v\|_{L^2(B)}.
\]
The equality $\ell_a(\one)=1$ follows from the definition of $m_a$.
\end{proof}

For completeness, the scattering coefficient in physical variables is
obtained as follows. The unitary rescaling
\[
 (\mathscr U_jv)(\xi)=a^{3/2}v(z_j+a\xi)
\]
satisfies $\mathscr U_jN_{D_j}\mathscr U_j^{-1}=a^2N_B$.
Let
\begin{equation}\label{eq:physical-W}
 W_j=(\kappa_a^{-1}I-N_{D_j})^{-1}\one.
\end{equation}
Since $(\kappa_a^{-1}I-N_{D_j})^{-1}
=\kappa_a(I-\kappa_aN_{D_j})^{-1}$, the rescaled function is
\[
 W_j(z_j+a\xi)=\theta a^{-2}\psi_a(\xi).
\]
We get then the properties of the scattering coefficient in the original variables:
\begin{equation}\label{eq:W-estimates}
 \int_{D_j}W_j=\alpha_a,
 \qquad \|W_j\|_{L^2(D_j)}
 =\theta a^{-1/2}\|\psi_a\|_{L^2(B)}
 \le Ct a^{-1/2-h}.
\end{equation}
The projection $\Pi_*$ in the proof includes the whole eigenspace.
Therefore no simplicity assumption on $\lambda_*$ is needed. Only its
coupling to the constant function, expressed by $b_*>0$, is used.
\section{Estimates for the distribution of the centers}
\label{sec:distribution}

Throughout Sections~\ref{sec:distribution}--\ref{sec:boundary}, we fix
$a$ and a realization $\varpi\in\cE_a$. All constants are independent
of this choice. We first derive consequences of separation and of
Assumption~\ref{ass:LLN}.

\subsection{Packing and quadrature estimates}
\begin{lemma}\label{lem:packing}
For $r>0$ and $x\in\R^3$,
\begin{equation}\label{eq:packing-count}
 \#\{j:|z_j-x|<r\}\le C(1+r/d)^3.
\end{equation}
For every fixed $s>0$,
\begin{equation}\label{eq:packing-sums}
 w\sum_{j\ne i}|z_i-z_j|^{-s}\le
 \begin{cases}
 C_s,&0<s<3,\\
 C(1+|\log d|),&s=3,\\
 C_s d^{3-s},&s>3.
 \end{cases}
\end{equation}
The same bounds hold for column sums. In addition, for $r\ge c d$,
\begin{equation}\label{eq:near-sum}
 w\sum_{0<|z_j-z_i|<r}|z_i-z_j|^{-1}\le Cr^2.
\end{equation}
\end{lemma}
\begin{proof}
The balls $B(z_j,c_0d/3)$ are disjoint. If $|z_j-x|<r$, each of these
balls is contained in $B(x,r+c_0d/3)$. Comparison of their volumes gives
\eqref{eq:packing-count}. Choose $c_1>0$ smaller than $c_0$ and divide
the distances from $z_i$ into dyadic shells
\[
 2^kc_1d\le |z_j-z_i|<2^{k+1}c_1d,
 \qquad k=0,\ldots,K,
\]
where $2^K d$ is comparable to $\operatorname{diam}\Omega$.
There are at most $C2^{3k}$ centers in the $k$th shell. Since $w=d^3$,
its contribution to the weighted sum in \eqref{eq:packing-sums} is at
most $C_s d^{3-s}2^{(3-s)k}$. Summing this geometric series gives the
three alternatives. Restricting the sum to shells with $2^kd\lesssim r$
and taking $s=1$ gives \eqref{eq:near-sum}. The column estimates are
identical, since the distance is symmetric. This dyadic-shell argument,
and closely related distribution-uniform interaction estimates, are also
used in \cite{BouzekriSini2021}.
\end{proof}

\begin{lemma}\label{lem:quadrature}
For $\psi\in C^1(\overline\Omega)$,
\begin{equation}\label{eq:green-quadrature}
 \max_i\left|
 w\sum_{j\ne i}G(z_i,z_j)\psi(z_j)
 -\int_\Omega G(z_i,y)\rho(y)\psi(y)\dd y\right|
 \le C(D_a^{1/2}+d^2)\|\psi\|_{C^1(\overline\Omega)}.
\end{equation}
\end{lemma}
\begin{proof}
Take a smooth cutoff $\chi_r$ which vanishes on $B(z_i,r)$ and equals
one outside $B(z_i,2r)$, with $|\nabla\chi_r|\le C/r$, where
$r\ge c d$. The estimates \eqref{eq:green-bounds} imply
\[
 \|\chi_r(\cdot)G(z_i,\cdot)\psi\|_{C^1(\overline\Omega)}
 \le Cr^{-2}\|\psi\|_{C^1(\overline\Omega)}
\]
for $r$ bounded above by a fixed constant. Assumption~\ref{ass:LLN}
therefore bounds the discrepancy of the cutoff expression by
$CD_ar^{-2}\|\psi\|_{C^1}$. Its omitted discrete part is bounded by
$Cr^2\|\psi\|_\infty$, by \eqref{eq:near-sum}. The omitted integral
has the same bound, since
\[
 \int_{\Omega\cap B(z_i,2r)}|z_i-y|^{-1}\dd y\le Cr^2
\]
and $\rho$ is bounded. The complete difference is consequently bounded
by $C(r^2+D_ar^{-2})\|\psi\|_{C^1}$. Taking
$r=\max\{c d,D_a^{1/4}\}$ proves \eqref{eq:green-quadrature}.
If $D_a=0$, take $r=c d$.
\end{proof}

The estimate \eqref{eq:green-quadrature} gives the expected density-weighted
integral equation for smooth functions. For a uniform inverse estimate,
we also need a comparison on Sobolev spaces. We establish it by replacing
each atom with a constant density on a small ball.

\subsection{An auxiliary absolutely continuous measure}
Choose $b_0>0$, depending only on $c_0,c_b$, so that the balls
$B(z_j,4b_0d)$ are disjoint and contained in $\Omega$. Set
\begin{equation}\label{eq:averaged-measure}
 B_j^d=B(z_j,b_0d),\quad v_d=|B_j^d|,\quad
 \rho_a^{\av}=\frac{w}{v_d}\sum_j\one_{B_j^d},\quad
 \nu_a=\rho_a^{\av}\dd x.
\end{equation}
Thus $\nu_a(B_j^d)=w$ and $\|\rho_a^{\av}\|_\infty\le C$.
For every $\varphi\in C^1(\overline\Omega)$,
\begin{align}
 \left|\int_\Omega\varphi\dd\nu_a
           -\int_\Omega\rho\varphi\dd x\right|
 &\le\left|w\sum_j\frac1{v_d}
       \int_{B_j^d}(\varphi(x)-\varphi(z_j))\dd x\right|
       +\mathfrak D_a\|\varphi\|_{C^1} \notag\\
 &\le C(D_a+d)\|\varphi\|_{C^1}.
 \label{eq:averaged-discrepancy}
\end{align}
Here we used $wM\le C_M$. These balls are used only to compare operators;
they are not substituted for the droplets in the acoustic problem.

\begin{lemma}\label{lem:form-discrepancy}
With $\epsilon_a=(D_a+d)^{1/5}$, one has
\begin{equation}\label{eq:form-discrepancy}
 \left|\int_\Omega(\rho_a^{\av}-\rho)u\overline v\dd x\right|
 \le C\epsilon_a
 \bigl(\|u\|_{H^1}\|v\|_{L^2}
           +\|u\|_{L^2}\|v\|_{H^1}\bigr),
 \qquad u,v\in H^1(\Omega).
\end{equation}
\end{lemma}
\begin{proof}
Set $\delta_a=D_a+d$. Since $\Omega$ is a bounded $C^{2,\alpha}$
domain, there is a linear extension operator
\[
 \mathcal E:H^1(\Omega)\longrightarrow H^1(\R^3)
\]
which is also bounded from $L^2(\Omega)$ to $L^2(\R^3)$. Let
$\eta\in C_c^\infty(B(0,1))$ satisfy $\eta\ge0$ and
$\int\eta=1$, and write $\eta_r(x)=r^{-3}\eta(x/r)$. For
$0<r<r_1$ define
\[
 u_r=\eta_r*(\mathcal Eu),\qquad
 v_r=\eta_r*(\mathcal Ev).
\]
We use the same symbols for their restrictions to $\Omega$.

We first record the smoothing estimates with their scaling. Since
$\int\eta_r=1$, translation estimates in $H^1(\R^3)$ give
\begin{align}
 \|u-u_r\|_{L^2(\Omega)}
 &\le Cr\|u\|_{H^1(\Omega)},
 &\|u_r\|_{L^2(\Omega)}
 &\le C\|u\|_{L^2(\Omega)}. \label{eq:mollifier-L2}
\end{align}
Young's inequality in the form
$\|f*g\|_\infty\le\|f\|_2\|g\|_2$ gives
\[
 \|u_r\|_\infty
 \le \|\eta_r\|_2\|\mathcal Eu\|_2
 \le Cr^{-3/2}\|u\|_2,
\]
while
$\|\nabla\eta_r\|_2=r^{-5/2}\|\nabla\eta\|_2$ yields
\[
 \|\nabla u_r\|_\infty
 \le Cr^{-5/2}\|u\|_2.
\]
The same estimates hold for $v_r$. Therefore
\begin{align}
 \|u_r\overline{v_r}\|_\infty
 &\le Cr^{-3}\|u\|_2\|v\|_2,\notag\\
 \|\nabla(u_r\overline{v_r})\|_\infty
 &\le \|\nabla u_r\|_\infty\|v_r\|_\infty
      +\|u_r\|_\infty\|\nabla v_r\|_\infty\notag\\
 &\le Cr^{-4}\|u\|_2\|v\|_2.\label{eq:mollified-gradient}
\end{align}
After decreasing $r_1$ so that $r\le1$, the $r^{-4}$ term dominates,
and hence
\begin{equation}\label{eq:mollified-product}
 \|u_r\overline{v_r}\|_{C^1(\overline\Omega)}
 \le Cr^{-4}\|u\|_{L^2(\Omega)}\|v\|_{L^2(\Omega)}.
\end{equation}

We now decompose the desired form. Put
$q_a=\rho_a^{\av}-\rho$. From
\eqref{eq:averaged-measure} and the bounds on $\rho$,
$\|q_a\|_\infty\le C$ uniformly in $a$. Thus
\begin{align*}
 \left|\int_\Omega q_a
   (u\overline v-u_r\overline{v_r})\dd x\right|
 &\le C\|u-u_r\|_2\|v\|_2
      +C\|u_r\|_2\|v-v_r\|_2\\
 &\le Cr\bigl(\|u\|_{H^1}\|v\|_2
          +\|u\|_2\|v\|_{H^1}\bigr),
\end{align*}
where \eqref{eq:mollifier-L2} was used in the last line.

For the smooth product, note that
\[
 \int_\Omega q_a u_r\overline{v_r}\dd x
 =\int_\Omega u_r\overline{v_r}\dd\nu_a
  -\int_\Omega\rho u_r\overline{v_r}\dd x.
\]
Hence \eqref{eq:averaged-discrepancy} and
\eqref{eq:mollified-product} imply
\[
 \left|\int_\Omega q_a u_r\overline{v_r}\dd x\right|
 \le C\delta_a r^{-4}\|u\|_2\|v\|_2.
\]
Combining the last two estimates gives, for every $0<r<r_1$,
\begin{align}\label{eq:form-discrepancy-preopt}
 \left|\int_\Omega q_a u\overline v\dd x\right|
 &\le Cr\bigl(\|u\|_{H^1}\|v\|_2
          +\|u\|_2\|v\|_{H^1}\bigr)
    +C\delta_a r^{-4}\|u\|_2\|v\|_2.
\end{align}
For small $a$, choose $r=\delta_a^{1/5}<r_1$. Then
$\delta_a r^{-4}=r=\delta_a^{1/5}=\epsilon_a$, and
$\|u\|_2\|v\|_2$ is bounded by either mixed norm appearing on the
right-hand side of \eqref{eq:form-discrepancy}. Thus
\eqref{eq:form-discrepancy-preopt} gives
\eqref{eq:form-discrepancy} after decreasing the asymptotic threshold
$a_0$ if necessary.
\end{proof}

\subsection{Coercivity and comparison of the auxiliary fields}
For $P\ge1$, use the norm
\begin{equation}\label{eq:P-norm}
 \|u\|_P^2=\|\nabla u\|_{L^2(\Omega)}^2
                      +P^2\|u\|_{L^2(\Omega)}^2.
\end{equation}
Let $q_* =\sup_{\omega\in I}\|q_\omega\|_\infty$. The two forms under
consideration are
\begin{align*}
 \mathfrak a_\rho(u,v)
   &=\int_\Omega\nabla u\cdot\nabla\overline v
                  +(P^2\rho-q_\omega)u\overline v,\\
 \mathfrak a_a(u,v)
   &=\int_\Omega\nabla u\cdot\nabla\overline v
                  +(P^2\rho_a^{\av}-q_\omega)u\overline v.
\end{align*}

\begin{lemma}\label{lem:coercivity}
There are $P_0\ge1$ and $c>0$ such that, if $P\ge P_0$ and
$P\epsilon_a\le c_{\mathrm{adm}}$ with $c_{\mathrm{adm}}$ sufficiently
small, then
\begin{equation}\label{eq:coercivity}
 \mathfrak a_\rho(u,u)\ge c\|u\|_P^2,
 \qquad \mathfrak a_a(u,u)\ge c\|u\|_P^2.
\end{equation}
Let $u_a^{\av,f}$ solve
\begin{equation}\label{eq:averaged-PDE}
 (-\Delta+P^2\rho_a^{\av}-q_\omega)u_a^{\av,f}=0,
 \qquad \partial_\nu u_a^{\av,f}=f,
 \qquad \Lambda_a^{\av}f=\tr u_a^{\av,f}.
\end{equation}
Then both this problem and \eqref{eq:effective-PDE} are uniquely
solvable, and
\begin{align}
 \|u_a^{\av,f}\|_P+\|u_{\rho,P}^f\|_P&\le C\|f\|_{\cX},
 \label{eq:averaged-energy}\\
 \|\Lambda_a^{\av}-\Lambda_{\rho,P}^{\eff}\|_{\cL(\cX,\cY)}
 &\le CP\epsilon_a.
 \label{eq:averaged-boundary-error}
\end{align}
\end{lemma}
\begin{proof}
Take $P_0^2\rho_*\ge2q_*$. Then
\[
 \mathfrak a_\rho(u,u)
 \ge\|\nabla u\|_2^2+\tfrac12\rho_*P^2\|u\|_2^2
 \ge c_1\|u\|_P^2,
 \qquad c_1=\min\{1,\rho_*/2\}.
\]
Lemma~\ref{lem:form-discrepancy}, together with
$\|u\|_{H^1}\le C\|u\|_P$ and
$\|u\|_2\le P^{-1}\|u\|_P$, yields
\begin{equation}\label{eq:form-P-bound}
 |\mathfrak a_a(u,v)-\mathfrak a_\rho(u,v)|
 \le CP\epsilon_a\|u\|_P\|v\|_P.
\end{equation}
Choosing $Cc_{\mathrm{adm}}\le c_1/2$ proves the second inequality
in \eqref{eq:coercivity}.

The trace theorem gives
$|\pair{f}{\tr v}|\le C\|f\|_{\cX}\|v\|_P$.
Lax--Milgram therefore gives both solutions and
\eqref{eq:averaged-energy}. Subtract their weak equations and put
$e=u_a^{\av,f}-u_{\rho,P}^f$. For $v\in H^1(\Omega)$,
\[
 \mathfrak a_a(e,v)
 =P^2\int_\Omega(\rho-\rho_a^{\av})
                        u_{\rho,P}^f\overline v.
\]
By \eqref{eq:form-P-bound}, the right-hand side is at most
$CP\epsilon_a\|u_{\rho,P}^f\|_P\|v\|_P$ in modulus. Take $v=e$,
use coercivity and \eqref{eq:averaged-energy}, and then apply the trace
theorem. This gives \eqref{eq:averaged-boundary-error}.
\end{proof}
\section{Invertibility of the algebraic system}
\label{sec:stability}

We prove \eqref{eq:uniform-FL-inverse} by comparing the matrix with an
integral operator associated with \eqref{eq:averaged-PDE}. This avoids a
smallness assumption on the full interaction matrix.

\subsection{The continuous inverse}
Let
\begin{equation}\label{eq:T-definition}
 \cH_a=L^2(\Omega,\nu_a),\qquad
 (T_aF)(x)=\int_\Omega G(x,y)F(y)\dd\nu_a(y).
\end{equation}
The space $\cH_a$ identifies functions which agree on the union of the
balls $B_j^d$. Since $\rho_a^{\av}$ is bounded and
$\sup_x\int_\Omega|x-y|^{-2}\dd y<\infty$, the kernel belongs to
$L^2(\nu_a\otimes\nu_a)$. Hence $T_a$ is bounded, compact and
self-adjoint on $\cH_a$.

\begin{lemma}\label{lem:continuous-inverse}
Under \eqref{eq:joint-regime}, after fixing $P_0$ and
$c_{\mathrm{adm}}$ as in Lemma~\ref{lem:coercivity},
\begin{equation}\label{eq:continuous-inverse}
 \|(I+tT_a)^{-1}\|_{\cH_a\to\cH_a}\le C,
 \qquad t=P^2.
\end{equation}
\end{lemma}
\begin{proof}
For $F\in\cH_a$, solve the zero-Neumann problem
\begin{equation}\label{eq:v-resolvent}
 (A_{0,\omega}+t\rho_a^{\av})v=\rho_a^{\av}F.
\end{equation}
The right-hand side is in $L^2(\Omega)$, because
\[
 \|\rho_a^{\av}F\|_{L^2(\Omega)}^2
 \le\|\rho_a^{\av}\|_\infty\|F\|_{\cH_a}^2.
\]
By coercivity and Cauchy--Schwarz,
\[
 c\|v\|_P^2
 \le\left|\int F\overline v\dd\nu_a\right|
 \le\|F\|_{\cH_a}\|v\|_{\cH_a}
 \le CP^{-1}\|F\|_{\cH_a}\|v\|_P.
\]
Consequently
\begin{equation}\label{eq:v-weight-bound}
 \|v\|_P\le CP^{-1}\|F\|_{\cH_a},
 \qquad t\|v\|_{\cH_a}\le C\|F\|_{\cH_a}.
\end{equation}
Applying $A_{0,\omega}^{-1}$ to \eqref{eq:v-resolvent} gives, on the
support of $\nu_a$,
\[
 v=T_a(F-tv),\qquad (I+tT_a)v=T_aF.
\]
It follows that $U=F-tv$ satisfies
\[
 (I+tT_a)U
 =F+tT_aF-t(v+tT_av)=F.
\]
The estimate \eqref{eq:v-weight-bound} yields
$\|U\|_{\cH_a}\le C\|F\|_{\cH_a}$.

To verify uniqueness, suppose $U+tT_aU=0$. Define the global function
$v=A_{0,\omega}^{-1}(\rho_a^{\av}U)$. On the support of $\nu_a$,
$v=T_aU$ and $U=-tv$. Thus
$(A_{0,\omega}+t\rho_a^{\av})v=0$ with zero Neumann data.
Coercivity gives $v=0$ and then $U=0$. This proves
\eqref{eq:continuous-inverse}.
\end{proof}

\subsection{Comparison with the point-interaction operator}
Define the constant extension and averaging operators by
\begin{equation}\label{eq:JQ}
 (JY)|_{B_i^d}=Y_i,\qquad
 (QF)_i=\frac1w\int_{B_i^d}F\dd\nu_a
       =\frac1{v_d}\int_{B_i^d}F\dd x.
\end{equation}
Then $J:\ell_w^2\to\cH_a$ is an isometry, $Q=J^*$, and $QJ=I$.
In particular, $JQ$ is the orthogonal projection onto the functions
which are constant on each ball. Put
\[
 (\mathsf K_aY)_i=w\sum_{j\ne i}G(z_i,z_j)Y_j,
 \qquad K_a=J\mathsf K_aQ.
\]
As an operator on $\cH_a$, $K_a$ has kernel $G(z_i,z_j)$ on
$B_i^d\times B_j^d$ when $i\ne j$, and zero on the diagonal blocks.

\begin{lemma}\label{lem:kernel-comparison}
One has
\begin{equation}\label{eq:kernel-comparison}
 \|T_a-K_a\|_{\cH_a\to\cH_a}\le Cd.
\end{equation}
Consequently \eqref{eq:uniform-FL-inverse} holds.
\end{lemma}
\begin{proof}
For $x\in B_i^d$, the integral of the absolute value of the kernel on
the diagonal block is bounded by
\[
 \int_{B_i^d}|G(x,y)|\dd\nu_a(y)
 \le C\int_{B(x,2b_0d)}|x-y|^{-1}\dd y\le Cd^2.
\]
If $x\in B_i^d$, $y\in B_j^d$ and $i\ne j$, our choice of $b_0$
ensures that every distance arising while moving $x$ to $z_i$ and $y$
to $z_j$ is comparable to $|z_i-z_j|$. The line segments lie inside
the corresponding balls. The gradient estimate
\eqref{eq:green-bounds} therefore gives
\[
 |G(x,y)-G(z_i,z_j)|\le Cd|z_i-z_j|^{-2}.
\]
Integration in $y$ and \eqref{eq:packing-sums}, with $s=2$, yield
\[
 \sup_x\int|G(x,y)-K_a(x,y)|\dd\nu_a(y)
 \le C(d^2+d)\le Cd.
\]
The same bound holds after interchanging $x$ and $y$. We spell out the
Schur argument. Let $L_a(x,y)$ be the kernel of $T_a-K_a$ with respect
to $\nu_a$. The preceding estimates say that
\begin{equation}\label{eq:schur-row-column}
 A_a:=\sup_x\int_\Omega|L_a(x,y)|\dd\nu_a(y)\le Cd,
 \qquad
 B_a:=\sup_y\int_\Omega|L_a(x,y)|\dd\nu_a(x)\le Cd.
\end{equation}
For completeness, the off-diagonal part of the first inequality is
\[
 \sum_{j\ne i}\int_{B_j^d}
    |G(x,y)-G(z_i,z_j)|\dd\nu_a(y)
 \le Cd\,w\sum_{j\ne i}|z_i-z_j|^{-2}\le Cd,
\]
for $x\in B_i^d$, while the diagonal block contributes only $Cd^2$.
The second inequality follows in the same way (or from the symmetry of
$G$ and of the block kernel).

If $F\in\cH_a$, Cauchy--Schwarz with the measure
$|L_a(x,y)|\dd\nu_a(y)$ gives, for a.e. $x$,
\[
 |(T_a-K_a)F(x)|^2
 \le\left(\int|L_a(x,y)|\dd\nu_a(y)\right)
     \left(\int|L_a(x,y)|\,|F(y)|^2\dd\nu_a(y)\right).
\]
Integrating in $x$, using Tonelli's theorem and
\eqref{eq:schur-row-column}, we obtain
\begin{align*}
 \|(T_a-K_a)F\|_{\cH_a}^2
 &\le A_a\int_\Omega |F(y)|^2
      \left(\int_\Omega|L_a(x,y)|\dd\nu_a(x)\right)\dd\nu_a(y)\\
 &\le A_aB_a\|F\|_{\cH_a}^2
 \le C d^2\|F\|_{\cH_a}^2.
\end{align*}
Thus $\|T_a-K_a\|\le\sqrt{A_aB_a}\le Cd$, which is
\eqref{eq:kernel-comparison}.

Now factor
\begin{equation}\label{eq:inverse-perturbation}
 I+tK_a=(I+tT_a)
 \bigl[I+t(I+tT_a)^{-1}(K_a-T_a)\bigr].
\end{equation}
The norm of the perturbation in brackets is at most $Ctd$.
For $d\le1$, we have
\begin{equation}\label{eq:small-parameters}
 d\le\epsilon_a^2,\qquad
 td\le(P\epsilon_a)^2\le c_{\mathrm{adm}}^2.
\end{equation}
Indeed $\epsilon_a^2=(D_a+d)^{2/5}\ge d^{2/5}\ge d$.
Choose $c_{\mathrm{adm}}$ so that $Ctd\le1/2$. A Neumann series in
\eqref{eq:inverse-perturbation}, together with
Lemma~\ref{lem:continuous-inverse}, gives a uniform inverse of
$I+tK_a$.

It remains to pass from the inverse on $\cH_a$ to the inverse of the
finite-dimensional Foldy--Lax matrix. Set
\[
 \mathcal V_a:=\operatorname{Ran}J,
 \qquad \Pi_a:=JQ.
\]
Since $Q=J^*$ and $QJ=I$, $\Pi_a$ is the orthogonal projection of
$\cH_a$ onto $\mathcal V_a$. Moreover,
\[
 K_a=J\mathsf K_aQ=\Pi_aK_a\Pi_a.
\]
In particular,
\begin{equation}\label{eq:Ka-range-properties}
 K_a(\mathcal V_a)\subset\mathcal V_a,
 \qquad K_a|_{\mathcal V_a^\perp}=0,
 \qquad (I-\Pi_a)K_a=K_a(I-\Pi_a)=0.
\end{equation}
Thus, with respect to the orthogonal decomposition
$\cH_a=\mathcal V_a\oplus\mathcal V_a^\perp$, the operator has the block
form
\[
 I+tK_a=
 \begin{pmatrix}
  (I+tK_a)|_{\mathcal V_a}&0\\
  0&I
 \end{pmatrix}.
\]
This already shows that its inverse preserves $\mathcal V_a$. Equivalently,
if $(I+tK_a)U\in\mathcal V_a$, then applying $I-\Pi_a$ and using
\eqref{eq:Ka-range-properties} gives
\[
 0=(I-\Pi_a)(I+tK_a)U=(I-\Pi_a)U,
\]
so $U\in\mathcal V_a$.

For $Y\in\ell_w^2$, the definitions give the exact intertwining identity
\begin{equation}\label{eq:J-intertwining}
 (I+tK_a)JY
 =JY+tJ\mathsf K_aQJY
 =J(I+t\mathsf K_a)Y
 =JB_{a,P}Y.
\end{equation}
Because $J$ is an isometric isomorphism from $\ell_w^2$ onto
$\mathcal V_a$, invertibility of $(I+tK_a)|_{\mathcal V_a}$ is equivalent
to invertibility of $B_{a,P}$, and
\[
 JB_{a,P}^{-1}
 =(I+tK_a)^{-1}J.
\]
Consequently, for every $Y\in\ell_w^2$,
\[
 \|B_{a,P}^{-1}Y\|_{\ell_w^2}
 =\|(I+tK_a)^{-1}JY\|_{\cH_a}
 \le C\|JY\|_{\cH_a}
 =C\|Y\|_{\ell_w^2}.
\]
This proves \eqref{eq:uniform-FL-inverse} with exactly the weighted norm
used in the statement.
\end{proof}

\subsection{Sampling the background solutions}
We shall use point values only for solutions of the background equation,
not for arbitrary $H^1$ functions.

\begin{lemma}\label{lem:sampling}
Let $s_i^f=S^f(z_i)$ and
$\widehat S_i^f(\xi)=S^f(z_i+a\xi)$ for $\xi\in B$.
For small $a$,
\begin{align}
 \|s^f\|_{\ell_w^2}
 +\left(w\sum_i\sup_{B_i^d}|\nabla S^f|^2\right)^{1/2}
 &\le C\|f\|_{\cX},\label{eq:sampling}\\
 \|S^f-Js^f\|_{\cH_a}&\le Cd\|f\|_{\cX},
 \label{eq:sampling-average}\\
 \left(w\sum_i\|\widehat S_i^f-s_i^f\|_{L^2(B)}^2\right)^{1/2}
 &\le Ca\|f\|_{\cX}.
 \label{eq:sampling-droplet}
\end{align}
\end{lemma}
\begin{proof}
Write $r=b_0d$. Interior estimates for $\Delta S^f+q_\omega S^f=0$
on $B(z_i,4r)$ give
\begin{align}
 |S^f(z_i)|&\le Cr^{-3/2}\|S^f\|_{L^2(B(z_i,2r))},
 \label{eq:local-value}\\
 \sup_{B(z_i,r)}|\nabla S^f|
 &\le Cr^{-3/2}\|S^f\|_{H^1(B(z_i,4r))}.
 \label{eq:local-gradient}
\end{align}
We recall why no additional inverse power of $r$ is needed in the
second estimate when the full $H^1$ norm is used. Differentiate the
equation in the weak sense:
\[
 (\Delta+q_\omega)\partial_kS^f
 =-(\partial_kq_\omega)S^f.
\]
The scaled interior $L^2$-to-$L^\infty$ estimate for this equation,
followed by the value estimate for $S^f$ on a larger ball, yields
\[
 \sup_{B(z_i,r)}|\partial_kS^f|
 \le C r^{-3/2}\|\partial_kS^f\|_{L^2(B(z_i,2r))}
       +Cr^2\|\nabla q_\omega\|_\infty
                    \sup_{B(z_i,2r)}|S^f|
 \le Cr^{-3/2}\|S^f\|_{H^1(B(z_i,4r))}.
\]
The constants are uniform since $r\le r_1$ and $q_\omega$ has a common
$W^{1,\infty}$ bound. These estimates can first be proved for smooth
coefficients and then passed to the limit; alternatively, interior
$W^{2,p}$ estimates give the same conclusion.

The balls $B(z_i,4r)$ are disjoint, and $wr^{-3}$ is a fixed constant.
Square \eqref{eq:local-value} and \eqref{eq:local-gradient}, multiply by
$w$, and sum in $i$. The background bound
$\|S^f\|_{H^1(\Omega)}\le C\|f\|_{\cX}$ proves
\eqref{eq:sampling}. On $B_i^d$,
\[
 |S^f(x)-S^f(z_i)|\le r\sup_{B_i^d}|\nabla S^f|.
\]
Integrating against $\nu_a$ and summing proves
\eqref{eq:sampling-average}. For small $a$, the set $z_i+aB$ and all
segments from $z_i$ to its points lie in $B_i^d$. The same argument,
with $|a\xi|\le a\sup_B|\xi|$, proves
\eqref{eq:sampling-droplet}. Convexity of the droplet is not required.
\end{proof}
\section{Derivation of the Point-Interaction approximation}
\label{sec:FL}

We next relate the exact field to the algebraic system of the preceding
section. The normalization of the unknowns is important: it removes the
local resonant factor from the remainder in the algebraic equation.

\subsection{The local integral equations}
By \eqref{eq:green-sign}, the exact equation \eqref{eq:exact-PDE} is
equivalent to
\begin{equation}\label{eq:exact-LS}
 u_a^f(x)=S^f(x)+\int_DG(x,y)(\kappa_a-q_\omega(y))u_a^f(y)\dd y,
 \qquad x\in\Omega.
\end{equation}
For the moment, suppose that a solution exists. We shall prove an a
priori estimate and then obtain existence by Fredholm theory. Write
\[
 \widehat u_i(\xi)=u_a^f(z_i+a\xi),\qquad
 \|\widehat u\|_{\mathbb H_a}^2
   =w\sum_i\|\widehat u_i\|_{L^2(B)}^2,
 \qquad U=\|\widehat u\|_{\mathbb H_a}.
\]
The physical and rescaled norms are related by
\begin{equation}\label{eq:physical-rescaled-norm}
 \|u_a^f\|_{L^2(D)}=a^{3/2}w^{-1/2}U.
\end{equation}
Define the matrix $\mathsf G^H$ by
\[
 \mathsf G^H_{ij}=
 \begin{cases}
 G(z_i,z_j),&i\ne j,\\
 H(z_i,z_i),&i=j.
 \end{cases}
\]
Restrict \eqref{eq:exact-LS} to $D_i$, use $G=\Phi_0+H$ on its
self-interaction term, and change variables to $B$. This gives the
exact identity
\begin{equation}\label{eq:local-LS}
 (I-\theta N_B)\widehat u_i
 =\one\left[s_i^f+a\theta\sum_j\mathsf G^H_{ij}
                            \int_B\widehat u_j\right]+F_i,
\end{equation}
where $F_i=F_i^{\mathrm{inc}}+F_i^{\mathrm{off}}
+F_i^{\mathrm{self}}+F_i^{\mathrm{bg}}$, with
\begin{align}
 F_i^{\mathrm{inc}}(\xi)
 &=\widehat S_i^f(\xi)-s_i^f,\label{eq:Finc}\\
 F_i^{\mathrm{off}}(\xi)
 &=a\theta\sum_{j\ne i}\int_B
   [G(z_i+a\xi,z_j+a\eta)-G(z_i,z_j)]\widehat u_j(\eta)\dd\eta,
 \label{eq:Foff}\\
 F_i^{\mathrm{self}}(\xi)
 &=a\theta\int_B
   [H(z_i+a\xi,z_i+a\eta)-H(z_i,z_i)]\widehat u_i(\eta)\dd\eta,
 \label{eq:Fself}\\
 F_i^{\mathrm{bg}}(\xi)
 &=-a^3\sum_j\int_B G(z_i+a\xi,z_j+a\eta)
        q_\omega(z_j+a\eta)\widehat u_j(\eta)\dd\eta.
 \label{eq:Fbg}
\end{align}
No interaction between distinct droplets has been discarded in
\eqref{eq:local-LS}.

\begin{lemma}\label{lem:local-remainders}
The remainder satisfies
\begin{equation}\label{eq:F-bound}
 \|F\|_{\mathbb H_a}\le Ca\|f\|_{\cX}+Ca^{1+h}U.
\end{equation}
\end{lemma}
\begin{proof}
The incident part is bounded by \eqref{eq:sampling-droplet}. Put
$r_{ij}=|z_i-z_j|$. Since $a/d\to0$, all the off-diagonal distances
in \eqref{eq:Foff} are comparable to $r_{ij}$. Taylor's formula and
\eqref{eq:green-bounds} give the block-operator bound
\[
 \|F_i^{\mathrm{off}}\|_{L^2(B)}
 \le Ca^2\sum_{j\ne i}r_{ij}^{-2}
                              \|\widehat u_j\|_{L^2(B)}.
\]
By \eqref{eq:packing-sums}, the row and column sums of this scalar
majorant are at most $Ca^2d^{-3}$. The block version of Schur's test
proved in Lemma~\ref{lem:block-schur} consequently gives
\[
 \|F^{\mathrm{off}}\|_{\mathbb H_a}\le Ca^2d^{-3}U.
\]
The Lipschitz estimate for $H$, with $b$ comparable to $d$, yields
\[
 \|F^{\mathrm{self}}\|_{\mathbb H_a}\le Ca^2d^{-2}U.
\]
For the off-diagonal part of \eqref{eq:Fbg}, the scalar majorant is
$Ca^3r_{ij}^{-1}$; its row and column sums are at most $Ca^3d^{-3}$.
On the diagonal, the Newtonian part is $a^2N_B$ followed by
multiplication by the bounded function $q_\omega(z_i+a\cdot)$.
The regular part has norm at most $Ca^3d^{-1}$. Hence
\[
 \|F^{\mathrm{bg}}\|_{\mathbb H_a}
 \le C(a^3d^{-3}+a^2+a^3d^{-1})U.
\]
For $a,d\le1$, each displayed coefficient is bounded by $Ca^2d^{-3}$.
Since $d^3=a^{1-h}$, we have $a^2d^{-3}=a^{1+h}$.
Adding the four estimates proves \eqref{eq:F-bound}.
\end{proof}

\subsection{Normalized amplitudes and the inverse estimate}
Set
\begin{equation}\label{eq:Y-definition}
 Y_i=\frac1{m_a}\int_B\widehat u_i
     =\frac{\kappa_a}{\alpha_a}\int_{D_i}u_a^f,
 \qquad r_i=\ell_a(F_i).
\end{equation}
The equality of the two expressions for $Y_i$ follows from
$\alpha_a=a\theta m_a$ and $\kappa_a=\theta a^{-2}$.
Apply $R_a$ to \eqref{eq:local-LS} and integrate on $B$. Division by
$m_a$ gives
\begin{equation}\label{eq:exact-FL}
 (I-\alpha_a\mathsf G^H)Y=s^f+r,
 \qquad \|r\|_{\ell_w^2}\le C\|F\|_{\mathbb H_a}.
\end{equation}
The uniform bound on $\ell_a$ in Lemma~\ref{lem:local} is used in the
last inequality.

\begin{lemma}\label{lem:physical-matrix}
Under \eqref{eq:joint-regime}, after reducing $a_0$ and
$c_{\mathrm{adm}}$ if necessary,
\begin{equation}\label{eq:physical-matrix-inverse}
 \|(I-\alpha_a\mathsf G^H)^{-1}\|_{\ell_w^2\to\ell_w^2}\le C.
\end{equation}
Moreover, with $\beta_i=1-\alpha_aH(z_i,z_i)$,
\begin{equation}\label{eq:beta-estimate}
 \max_i|\beta_i-1|\le Ct d^2.
\end{equation}
\end{lemma}
\begin{proof}
The off-diagonal difference between $I-\alpha_a\mathsf G^H$ and
$B_{a,P}$ is $-(\alpha_a+tw)G(z_i,z_j)$. Since
$|\alpha_a+tw|\le Ca$ and $\|wG_{\mathrm{off}}\|_{\ell_w^2\to\ell_w^2}
\le C$ by packing and Schur's test, this part has norm at most
$Ca/w=Ca^h$. On the diagonal, \eqref{eq:H-bounds} gives
\[
 |\alpha_aH(z_i,z_i)|\le Ctw d^{-1}=Ct d^2.
\]
This proves \eqref{eq:beta-estimate} and
\begin{equation}\label{eq:matrix-difference}
 \|I-\alpha_a\mathsf G^H-B_{a,P}\|
 \le C(a^h+t d^2).
\end{equation}
Combine this estimate with the uniform inverse of $B_{a,P}$.
The perturbation is smaller than $1/2$ for small $a$ and
$c_{\mathrm{adm}}$, using \eqref{eq:small-parameters}.
The resolvent identity, or a Neumann series, proves
\eqref{eq:physical-matrix-inverse}.
\end{proof}

\begin{proposition}\label{prop:physical-FL}
The exact problem \eqref{eq:exact-PDE} is uniquely solvable on $\cE_a$
in the regime \eqref{eq:joint-regime}. If
$Z^f=B_{a,P}^{-1}s^f$, then
\begin{align}
 U&\le Ct a^{-h}\|f\|_{\cX},&
 \|u_a^f\|_{L^2(D)}&\le Ct a^{1-h/2}\|f\|_{\cX},
 \label{eq:physical-field}\\
 \|F\|_{\mathbb H_a}+\|r\|_{\ell_w^2}
 &\le Ct a\|f\|_{\cX},&
 \|Y\|_{\ell_w^2}+\|Z^f\|_{\ell_w^2}&\le C\|f\|_{\cX},
 \label{eq:physical-amplitudes}\\
 \|Y-Z^f\|_{\ell_w^2}
 &\le C(ta+a^h+t d^2)\|f\|_{\cX}.
 \label{eq:FL-amplitude-error}
\end{align}
\end{proposition}
\begin{proof}
The scalar bracket in \eqref{eq:local-LS} equals $Y_i-\ell_a(F_i)$,
by \eqref{eq:exact-FL}. Therefore the exact reconstruction of the local
field is
\begin{equation}\label{eq:local-reconstruction}
 \widehat u_i=\psi_aY_i+R_aF_i-\psi_a\ell_a(F_i).
\end{equation}
The last two terms have zero integral on $B$. From
Lemmas~\ref{lem:local} and \ref{lem:physical-matrix}, and
\eqref{eq:sampling}, we obtain
\begin{align*}
 U&\le Ct a^{-h}\bigl(\|Y\|_{\ell_w^2}+\|F\|_{\mathbb H_a}\bigr),\\
 \|Y\|_{\ell_w^2}&\le C\bigl(\|f\|_{\cX}+\|F\|_{\mathbb H_a}\bigr).
\end{align*}
Substituting \eqref{eq:F-bound} into these inequalities gives
\begin{equation}\label{eq:field-absorption}
 U\le Ct a^{-h}\|f\|_{\cX}+Cta\,U.
\end{equation}
Since $a\le d$ for small $a$, the last coefficient is bounded by
$Ctd\le Cc_{\mathrm{adm}}^2$. Absorbing it proves the first estimate in
\eqref{eq:physical-field}. The second follows from
\eqref{eq:physical-rescaled-norm}, because
$a^{3/2}w^{-1/2}a^{-h}=a^{1-h/2}$.
Equation \eqref{eq:F-bound} now yields
$\|F\|_{\mathbb H_a}\le C(a+ta)\|f\|_{\cX}
\le Cta\|f\|_{\cX}$. This proves the estimate for $r$ as well.
Since $ta$ is small, the estimate for $Y$ follows from
\eqref{eq:exact-FL}; the one for $Z^f$ follows from
\eqref{eq:uniform-FL-inverse} and \eqref{eq:sampling}.

For the difference of the amplitudes, subtract
$B_{a,P}Z^f=s^f$ from \eqref{eq:exact-FL} to obtain
\[
 B_{a,P}(Y-Z^f)
 =r+[B_{a,P}-(I-\alpha_a\mathsf G^H)]Y.
\]
Use \eqref{eq:matrix-difference} and the estimates just proved to
obtain \eqref{eq:FL-amplitude-error}.

It remains to justify existence. For a fixed $a$ and a fixed
configuration, the integral operator in \eqref{eq:exact-LS}, restricted
to $L^2(D)$, is compact. Indeed its kernel is a bounded multiplier,
depending on $a$, times $G(x,y)$, and $|x-y|^{-2}$ is integrable on
$D\times D$. If $f=0$, the a priori inequality
\eqref{eq:field-absorption} gives $U=0$. Thus the homogeneous integral
equation has only the zero solution. Fredholm's alternative gives a
unique solution on $D$ for every $S^f|_D$.

Extend this solution to $\Omega$ by the right-hand side of
\eqref{eq:exact-LS}. The source
$(\kappa_a-q_\omega)\one_Du_a^f$ lies in $L^2(\Omega)$ for each fixed
$a$. The background resolvent maps it into $H^2(\Omega)$ with zero
Neumann trace. The extended field therefore belongs to $H^1(\Omega)$,
has flux $f$, agrees with the solution on $D$, and satisfies
\eqref{eq:exact-PDE} distributionally. Conversely any weak solution of
that PDE satisfies \eqref{eq:exact-LS} by the background inverse.
This proves existence and uniqueness, and completes the proof.
\end{proof}

\subsection{Closure of the point-interaction approximation}
We summarize explicitly how Lemmas~\ref{lem:local-remainders} and
\ref{lem:physical-matrix} and Proposition~\ref{prop:physical-FL} close
the argument of this section. Lemma~\ref{lem:local-remainders} reduces all
errors made in replacing the exact local integral equations by their
monopole parts to
\[
 \|F\|_{\mathbb H_a}\le Ca\|f\|_{\cX}+Ca^{1+h}U.
\]
Lemma~\ref{lem:physical-matrix} then supplies a uniform inverse for the
actual amplitude matrix $I-\alpha_a\mathsf G^H$. These two facts are the
inputs in the absorption argument of Proposition~\ref{prop:physical-FL},
which gives, without assuming solvability in advance,
\begin{equation}\label{eq:section5-closed-estimates}
 U\le Ct a^{-h}\|f\|_{\cX},\qquad
 \|r\|_{\ell_w^2}\le Cta\|f\|_{\cX},\qquad
 \|Y\|_{\ell_w^2}\le C\|f\|_{\cX}.
\end{equation}
The same proposition proves existence and uniqueness of the exact physical
field by Fredholm's alternative, so no additional solvability argument is
needed later.

The remaining step is to compare the physical amplitudes $Y$ with the
ideal Foldy--Lax amplitudes $Z^f=B_{a,P}^{-1}s^f$. Subtracting the two
systems gives the exact identity
\begin{equation}\label{eq:section5-amplitude-identity}
 B_{a,P}(Y-Z^f)
 =r+\bigl[B_{a,P}-(I-\alpha_a\mathsf G^H)\bigr]Y.
\end{equation}
Using \eqref{eq:uniform-FL-inverse},
\eqref{eq:matrix-difference}, and
\eqref{eq:section5-closed-estimates} yields
\[
 \|Y-Z^f\|_{\ell_w^2}
 \le C(ta+a^h+td^2)\|f\|_{\cX},
\]
which is \eqref{eq:FL-amplitude-error}. Thus the exact many-droplet problem
has been reduced, with a quantified error, to the ideal algebraic system
$B_{a,P}Z^f=s^f$.

It is useful to record already the form in which this estimate enters the
boundary-map comparison. For another boundary datum $g$, the sampling
estimate gives $\|s^g\|_{\ell_w^2}\le C\|g\|_{\cX}$. Since
$\alpha_a=-tw+O(a)$, weighted Cauchy--Schwarz gives
\begin{align}
 &\left|\alpha_a\sum_iY_i^fs_i^g
          +tw\sum_iZ_i^fs_i^g\right|\notag\\
 &\qquad\le
 t\|Y^f-Z^f\|_{\ell_w^2}\|s^g\|_{\ell_w^2}
 +\frac{|\alpha_a+tw|}{w}
       \|Y^f\|_{\ell_w^2}\|s^g\|_{\ell_w^2}\notag\\
 &\qquad\le
 C\bigl(t^2a+t a^h+t^2d^2\bigr)
   \|f\|_{\cX}\|g\|_{\cX}.\label{eq:section5-boundary-ready}
\end{align}
Here $|\alpha_a+tw|/w\le Ca/w=Ca^h$. Equation
\eqref{eq:section5-boundary-ready} is the point-interaction error needed
in Section~\ref{sec:boundary}; that section only has to add the local
variation of the test field and compare the ideal discrete system with
the auxiliary continuous medium.
\section{Proof of Theorem~\ref{thm:main}: Convergence of the boundary maps}
\label{sec:boundary}

We compare the exact boundary map first with a map defined by the
Foldy--Lax amplitudes, and then with the auxiliary and effective maps.
All pairings in this section are bilinear. They are well-defined for
$f,g\in\cX$ by the trace theorem; the identities may also be obtained
first for smooth data and then extended by continuity.

\subsection{The exact boundary identity}
The weak forms of the exact and background equations give
\begin{equation}\label{eq:boundary-identity}
 \pair{(\Lambda_{a,P}-\Lambda_{0,\omega})f}{g}
 =\int_D(\kappa_a-q_\omega)u_a^f S^g\dd x.
\end{equation}
Indeed, the exact equation implies
\[
 \mathfrak b_0(u_a^f,S^g)
 =\pair{f}{\tr S^g}+\int_D(\kappa_a-q_\omega)u_a^fS^g.
\]
On the other hand, the background equation for $S^g$ gives
$\mathfrak b_0(u_a^f,S^g)=\pair{g}{\tr u_a^f}$.
The background map is symmetric, because
$\mathfrak b_0(S^f,S^g)=\mathfrak b_0(S^g,S^f)$.
Subtracting its boundary pairing proves
\eqref{eq:boundary-identity}.

Using \eqref{eq:Y-definition}, expand the right-hand side as
\begin{equation}\label{eq:boundary-discrete}
 \pair{(\Lambda_{a,P}-\Lambda_{0,\omega})f}{g}
 =\alpha_a\sum_iY_i^f s_i^g+\mathcal R_a(f,g).
\end{equation}
We claim that
\begin{equation}\label{eq:boundary-remainder}
 |\mathcal R_a(f,g)|\le Cta\|f\|_{\cX}\|g\|_{\cX}.
\end{equation}
Indeed, the part containing the variation of $S^g$ on the droplets is
bounded, after rescaling, by
\begin{align*}
 \kappa_a a^3\sum_i
  \|\widehat u_i^f\|_{L^2(B)}
  \|\widehat S_i^g-s_i^g\|_{L^2(B)}
 &\le C\frac a w U
   \left(w\sum_i\|\widehat S_i^g-s_i^g\|_2^2\right)^{1/2}\\
 &\le C\frac{a^2}w\,t a^{-h}
                           \|f\|_{\cX}\|g\|_{\cX}
 =Cta\|f\|_{\cX}\|g\|_{\cX}.
\end{align*}
Here we used \eqref{eq:sampling-droplet} and
\eqref{eq:physical-field}. The part containing $q_\omega$ is at most
\[
 C\frac{a^3}w U
 \left(w\sum_i\|\widehat S_i^g\|_2^2\right)^{1/2}
 \le Ct a^2\|f\|_{\cX}\|g\|_{\cX},
\]
since \eqref{eq:sampling} and \eqref{eq:sampling-droplet} bound the
last factor. This proves \eqref{eq:boundary-remainder}.

\subsection{The Foldy--Lax boundary map}
Define $\Lambda_{a,P}^{\fl}:\cX\to\cY$ by
\begin{equation}\label{eq:FL-boundary-map}
 \pair{(\Lambda_{a,P}^{\fl}-\Lambda_{0,\omega})f}{g}
 =-tw\sum_iZ_i^f s_i^g,
 \qquad B_{a,P}Z^f=s^f.
\end{equation}
The right-hand side is a continuous linear functional of $g$, because
$\|Z^f\|_{\ell_w^2}\le C\|f\|_{\cX}$ and
$\|s^g\|_{\ell_w^2}\le C\|g\|_{\cX}$. It therefore defines an element
of $\cY$, with the usual identification of the dual of $\cX$.
From \eqref{eq:boundary-discrete},
\begin{align*}
 \pair{(\Lambda_{a,P}-\Lambda_{a,P}^{\fl})f}{g}
 &=-tw\sum_i(Y_i^f-Z_i^f)s_i^g
   +(\alpha_a+tw)\sum_iY_i^f s_i^g+\mathcal R_a(f,g).
\end{align*}
Cauchy--Schwarz in $\ell_w^2$, the coefficient expansion
\eqref{eq:alpha-uniform}, and
\eqref{eq:FL-amplitude-error} give
\begin{equation}\label{eq:exact-FL-boundary}
 \|\Lambda_{a,P}-\Lambda_{a,P}^{\fl}\|_{\cL(\cX,\cY)}
 \le C(t^2a+t a^h+t^2d^2).
\end{equation}
The factor in front of the coefficient remainder is
$|\alpha_a+tw|/w\le Ca^h$; it is bounded by $Ct a^h$ since $t\ge1$.
The term $Cta$ in \eqref{eq:boundary-remainder} is absorbed in $Ct^2a$.

\subsection{Comparison with the auxiliary map}
Let $U^f$ denote the restriction of $u_a^{\av,f}$ to the support of
$\nu_a$. Applying the background resolvent to
\eqref{eq:averaged-PDE} gives
\begin{equation}\label{eq:averaged-LS}
 (I+tT_a)U^f=S^f|_{\supp\nu_a},
 \qquad \|U^f\|_{\cH_a}\le C\|f\|_{\cX}.
\end{equation}
For the norm bound, one may use either
\eqref{eq:averaged-energy} or Lemma~\ref{lem:continuous-inverse} and
$\|S^f\|_{\cH_a}\le C\|S^f\|_{L^2(\Omega)}$.
Recall the map $J$ defined in (\ref{eq:JQ}). Since $(I+tK_a)JZ^f=Js^f$, subtraction yields
\[
 (I+tT_a)(U^f-JZ^f)
 =(S^f-Js^f)-t(T_a-K_a)JZ^f.
\]
Use \eqref{eq:sampling-average},
\eqref{eq:kernel-comparison}, the isometry of $J$, and the two inverse
estimates. We obtain
\begin{equation}\label{eq:averaged-FL-fields}
 \|U^f-JZ^f\|_{\cH_a}\le C(d+td)\|f\|_{\cX}
 \le Ctd\|f\|_{\cX}.
\end{equation}
The analogue of \eqref{eq:boundary-identity} for the auxiliary equation
is
\begin{equation}\label{eq:averaged-boundary-identity}
 \pair{(\Lambda_a^{\av}-\Lambda_{0,\omega})f}{g}
 =-t\int_\Omega U^fS^g\dd\nu_a.
\end{equation}
Also, $w\sum_iZ_i^fs_i^g=\int JZ^fJs^g\dd\nu_a$.
Subtract \eqref{eq:FL-boundary-map} from
\eqref{eq:averaged-boundary-identity} and use
\[
 U^fS^g-JZ^fJs^g
 =(U^f-JZ^f)S^g+JZ^f(S^g-Js^g).
\]
Cauchy--Schwarz, \eqref{eq:averaged-FL-fields} and
\eqref{eq:sampling-average} give
\begin{equation}\label{eq:averaged-FL-map}
 \|\Lambda_a^{\av}-\Lambda_{a,P}^{\fl}\|_{\cL(\cX,\cY)}
 \le Ct^2d.
\end{equation}

Combining \eqref{eq:exact-FL-boundary},
\eqref{eq:averaged-FL-map} and
\eqref{eq:averaged-boundary-error}, we find
\[
 \|\Lambda_{a,P}-\Lambda_{\rho,P}^{\eff}\|
 \le C(t^2a+t a^h+t^2d^2+t^2d+P\epsilon_a).
\]
Since $a\le d\le1$, this is \eqref{eq:main-error}. Existence and
uniqueness were proved in Proposition~\ref{prop:physical-FL}, and the
uniform matrix inverse was proved in Lemma~\ref{lem:kernel-comparison}.
This completes the proof of Theorem~\ref{thm:main}.

Finally, the limiting integral equation is
\begin{equation}\label{eq:effective-LS}
 u_{\rho,P}^f(x)+P^2\int_\Omega
       G_\omega(x,y)\rho(y)u_{\rho,P}^f(y)\dd y=S^f(x).
\end{equation}
Applying $A_{0,\omega}$ gives precisely
\eqref{eq:effective-PDE}. This also identifies the normalization of the
effective perturbation: the measure $w\sum_j\delta_{z_j}$ produces
$P^2\rho$, without an additional factor involving the total mass of
$\rho$.
\section{A class of random configurations}
\label{sec:random}

We now verify the geometric and distribution assumptions for a concrete
probability law. The construction permits random occupancies and random
displacements, but retains a deterministic lower bound on the separation.
Consequently, it is not the unrestricted sampling of independent centers from
$\rho(x)\dd x$. We call this the
\emph{Bernoulli-rounded lattice model}.

\subsection{Occupancy and displacement laws}
Let $\mathcal Q_d$ be the collection of cubes of side length $d$ in a
fixed cubic grid which are contained in $\Omega$ and whose distance
from $\partial\Omega$ is at least $2d$. Put
$\Omega_d=\bigcup_{Q\in\mathcal Q_d}Q$. The smoothness of the boundary
gives
\begin{equation}\label{eq:boundary-layer-volume}
 |\Omega\setminus\Omega_d|\le Cd.
\end{equation}
Indeed, a point at distance greater than $(2+\sqrt3)d$ from the boundary
belongs to an admissible cube. The complement is therefore contained
in a tubular neighborhood of $\partial\Omega$ of width $Cd$, whose
volume is at most $Cd$.

For $Q\in\mathcal Q_d$, set
\begin{equation}\label{eq:random-rounding}
 \lambda_Q=w^{-1}\int_Q\rho(x)\dd x,\qquad
 N_Q=\lfloor\lambda_Q\rfloor+B_Q,
 \qquad \PP(B_Q=1)=\lambda_Q-\lfloor\lambda_Q\rfloor.
\end{equation}
The Bernoulli variables $B_Q$ are independent. Since $w=d^3$,
\[
 \rho_*\le\lambda_Q\le\rho^*,\qquad
 \EE N_Q=\lambda_Q,
 \qquad 0\le N_Q\le K:=\lceil\rho^*\rceil+1.
\]
Fix $K$ points in the interior of the unit cube, with positive mutual
distance and positive distance from its boundary. In each cube $Q$,
scale and translate these points to obtain $K$ slots. Choose the radius
of the displacement balls to be a sufficiently small fixed multiple of
$d$. Conditional on $N_Q$, place one center in each of the first $N_Q$
slots, independently according to any prescribed law supported in its
displacement ball. Uniform distributions in these balls are one example.
The displacement variables can be defined for all $K$ slots in advance,
independently of the occupancy variables; unused slots are then ignored.
Here the occupancy law governs the number $N_Q$ of centers in
$Q$, whereas a displacement law specifies their random offsets from
fixed slots.

All centers in a given cube, and all centers in neighboring cubes,
are separated by at least $c_Kd$. They are also at distance at least
$c_bd$ from $\partial\Omega$. Moreover,
\[
 wM=w\sum_QN_Q\le Kd^3\#\mathcal Q_d\le C.
\]
The preceding assertions can be stated as an exact probability-one
statement. With the convention that a minimum over an empty set is
$+\infty$, define
\begin{align}
 \mathcal G_a:=\Bigl\{\varpi:\;&
 \min_{i\ne j}|z_i(\varpi)-z_j(\varpi)|\ge c_Kd,\quad
 \min_j\dist(z_j(\varpi),\partial\Omega)\ge c_bd,\notag\\
 & wM(a,\varpi)\le C_M\Bigr\}.\label{eq:geometric-event-random}
\end{align}
The slot construction implies $\mathcal G_a=\Xi$ (up to completion of
the probability space), and therefore
\begin{equation}\label{eq:geometry-probability-one}
 \PP(\mathcal G_a)=1.
\end{equation}
Thus Assumption~\ref{ass:hardcore} is satisfied with
$c_0=c_K$ and with probability exactly one, not only asymptotically one.
The probability space is the product of the occupancy and displacement
spaces over the cubes.
With $p_Q=\lambda_Q-\lfloor\lambda_Q\rfloor$ and
independent offset laws $\nu_{Q,j}$, its distribution at scale $a$ is
$\PP_a=\bigotimes_{Q\in\mathcal Q_d}
 (\operatorname{Ber}(p_Q)\otimes\bigotimes_{j=1}^K\nu_{Q,j})$.
The occupancy laws depend on $Q$ when $\rho$ is
nonconstant; no stationarity assumption is imposed. For a prescribed
countable sequence of scales, all these variables may be placed on one
product probability space.

The event used simultaneously in Assumptions~\ref{ass:hardcore} and
\ref{ass:LLN} will therefore be
\begin{equation}\label{eq:random-common-event-definition}
 \cE_a^{(m)}
 :=\mathcal G_a\cap
   \{\mathfrak D_a\le C_m d^{1/2}\}
 =\{\mathfrak D_a\le C_m d^{1/2}\}
 \quad\text{up to a null set}.
\end{equation}
Proposition~\ref{prop:random-discrepancy} below proves that, for small
$d$,
\[
 \PP(\cE_a^{(m)})\ge1-d^m.
\]
This is the precise probability entering \eqref{eq:common-event} for the
present random construction. Hence, one may take $p_a=d^m$ and
$D_a=C_m d^{1/2}$.

\subsection{A common-event estimate for the empirical measure}
\begin{proposition}\label{prop:random-discrepancy}
For every $m>0$, there are $C_m,d_m>0$ such that the preceding
construction satisfies
\begin{equation}\label{eq:random-discrepancy}
 \PP\{\mathfrak D_a\le C_m d^{1/2}\}\ge1-d^m,
 \qquad 0<d<d_m.
\end{equation}
The event is independent of $P$ (and of the boundary inputs).
\end{proposition}
\begin{proof}
Choose a point $\xi_Q\in Q$ and put $X_Q=N_Q-\lambda_Q$.
For $\|\varphi\|_{C^1(\overline\Omega)}\le1$, the displacement of
centers within $Q$ contributes at most $Cd\,wN_Q$. Also,
\[
 \left|\int_Q\rho(x)(\varphi(x)-\varphi(\xi_Q))\dd x\right|
 \le Cd\int_Q\rho(x)\dd x.
\]
Together with \eqref{eq:boundary-layer-volume}, these estimates give
\begin{equation}\label{eq:discrepancy-random-part}
 \left|\int\varphi\dd\mu_a-\int_\Omega\rho\varphi\dd x\right|
 \le Cd+\left|w\sum_{Q\in\mathcal Q_d}X_Q\varphi(\xi_Q)\right|.
\end{equation}
This estimate is uniform over the displacement variables.

To control the last term simultaneously for all $\varphi$, partition a
fixed box containing $\Omega$ into coarse cubes $R$ of side length
$L\in[d,1]$. Assign $Q$ to $R$ when $\xi_Q\in R$ and write
$n_R=\#\{Q:\xi_Q\in R\}$. The numbers of coarse cubes and fine cubes
satisfy
\begin{equation}\label{eq:coarse-counts}
 \#\{R:n_R>0\}\le CL^{-3},\qquad
 \sum_Rn_R\le Cd^{-3},\qquad n_R\le C(L/d)^3.
\end{equation}
Use a bounded $C^1$ extension of $\varphi$ to the fixed box. Such an
extension exists since $\Omega$ is smooth; its $C^1$ norm is at most
$C_\Omega$. Replace $\varphi(\xi_Q)$ by the value of this extension at
one point $\xi_R$ of the corresponding coarse cube. As $|X_Q|\le1$,
this replacement has error at most
$CLw\sum_Q1\le CL$. Consequently, from (\ref{eq:discrepancy-random-part}), we get
\begin{equation}\label{eq:uniform-discrepancy-bound}
 \mathfrak D_a\le C(d+L)+Cw\sum_R
                       \left|\sum_{Q:\xi_Q\in R}X_Q\right|.
\end{equation}
The right-hand side no longer depends on $\varphi$.

We now make the concentration step explicit. Write
\[
 p_Q:=\lambda_Q-\lfloor\lambda_Q\rfloor\in[0,1),
 \qquad X_Q=N_Q-\lambda_Q=B_Q-p_Q.
\]
Thus
\begin{equation}\label{eq:XQ-properties}
 \EE X_Q=0,\qquad |X_Q|\le1,\qquad
 \operatorname{Var}(X_Q)=p_Q(1-p_Q)\le\frac14,
\end{equation}
and the two possible values of $X_Q$ are $-p_Q$ and $1-p_Q$; hence the
length of its range is exactly one. The variables $X_Q$ remain
independent because the $B_Q$ are independent.

There are two equivalent ways to obtain the needed tail estimate. First,
Hoeffding's inequality for independent centered variables with range
length one \cite{Hoeffding} gives, for a fixed coarse cube $R$ with
$n_R>0$ and
$S_R:=\sum_{Q:\xi_Q\in R}X_Q$,
\begin{equation}\label{eq:hoeffding}
 \PP\{|S_R|>u\}
 \le2\exp\left(-\frac{2u^2}{n_R}\right)
 \le2\exp(-c u^2/n_R).
\end{equation}

For completeness, the same bound follows directly from exponential
moments. For every $s\in\R$,
\[
 \EE e^{sX_Q}
 =e^{-sp_Q}\bigl(1-p_Q+p_Qe^s\bigr)
 \le e^{s^2/8},
\]
where the last inequality is Hoeffding's elementary lemma for a centered
random variable supported in an interval of length one. Independence
then gives
\[
 \EE e^{sS_R}=\prod_{Q:\xi_Q\in R}\EE e^{sX_Q}
 \le e^{n_Rs^2/8}.
\]
For $s>0$, Markov's inequality yields
\[
 \PP\{S_R>u\}
 \le e^{-su}\EE e^{sS_R}
 \le\exp\left(-su+\frac{n_Rs^2}{8}\right).
\]
Minimizing the exponent at $s=4u/n_R$ gives
$\PP\{S_R>u\}\le e^{-2u^2/n_R}$. Applying the same computation to
$-S_R$ and adding the two tails gives \eqref{eq:hoeffding}.

Let $0<\delta<1$ and let
$N_L=\#\{R:n_R>0\}\le CL^{-3}$. Choose
\[
 u_R=C\sqrt{n_R\log(2N_L/\delta)}
\]
with $C$ large enough that \eqref{eq:hoeffding} gives
$\PP\{|S_R|>u_R\}\le\delta/N_L$. Summing over the at most $N_L$
occupied coarse cubes gives
\[
 \PP\left\{|S_R|\le u_R\ \text{for every occupied }R\right\}
 \ge1-\sum_{R:n_R>0}\frac{\delta}{N_L}
 \ge1-\delta.
\]
Since $N_L\le CL^{-3}$, on this event
\begin{equation}\label{eq:all-coarse-cubes}
 |S_R|\le C\sqrt{n_R\log(CL^{-3}/\delta)}
 \qquad\text{for every occupied coarse cube }R.
\end{equation}
Cauchy--Schwarz and \eqref{eq:coarse-counts} give
\[
 w\sum_R\sqrt{n_R}
 \le w\bigl(\#\{R:n_R>0\}\bigr)^{1/2}
           \left(\sum_Rn_R\right)^{1/2}
 \le C(d/L)^{3/2}.
\]
Inserting these inequalities into
\eqref{eq:uniform-discrepancy-bound}, we obtain, on the same event,
\begin{equation}\label{eq:coarse-discrepancy-rate}
 \mathfrak D_a
 \le C\left[d+L+(d/L)^{3/2}\sqrt{\log(CL^{-3}/\delta)}\right].
\end{equation}
Choose $L=d^{1/2}$ and $\delta=d^m$. The last term is bounded by
$C_m d^{3/4}\sqrt{1+|\log d|}\le C_m d^{1/2}$ for small $d$.
This proves \eqref{eq:random-discrepancy}.
\end{proof}

\begin{corollary}\label{cor:random-forward}
For each fixed $m>0$, the boundary maps for this construction satisfy,
with probability at least $1-d^m$,
\begin{equation}\label{eq:random-forward}
 \|\Lambda_{a,P}-\Lambda_{\rho,P}^{\eff}\|
 \le C_m\bigl(P^2a^h+P^4d+P d^{1/10}\bigr),
\end{equation}
simultaneously for all $P\ge P_0$ such that
$P d^{1/10}\le c_m$.
\end{corollary}
\begin{proof}
On the event in Proposition~\ref{prop:random-discrepancy}, take
$D_a=C_m d^{1/2}$ in Assumption~\ref{ass:LLN}. Then
$\epsilon_a=(D_a+d)^{1/5}\le C_m d^{1/10}$.
Theorem~\ref{thm:main} gives the assertion after adjusting $c_m$.
\end{proof}

We finally spell out the almost-sure statement. Fix $m>0$ and a
countable sequence of scales $d_k\downarrow0$, with the corresponding
$a_k=d_k^{3/(1-h)}$, and realize all occupancies and displacements on the
single product probability space described above. Define
\[
 E_k:=\{\mathfrak D_{a_k}\le C_m d_k^{1/2}\}.
\]
Proposition~\ref{prop:random-discrepancy} gives
$\PP(E_k^c)\le d_k^m$ for all sufficiently large $k$. Thus there is
$k_0$ such that, for the concrete choice $d_k=2^{-k}$,
\[
 \sum_{k=k_0}^\infty\PP(E_k^c)
 \le\sum_{k=k_0}^\infty2^{-mk}<\infty.
\]
Adding the finitely many terms with $k<k_0$ shows that
$\sum_{k=1}^\infty\PP(E_k^c)<\infty$. Hence the first
Borel--Cantelli lemma yields
\[
 \PP(E_k^c\ \text{infinitely often})=0,
\]
or, equivalently,
\begin{equation}\label{eq:random-typical-set}
 \Xi_{\rm typ}
 :=\bigcup_{K=1}^\infty\bigcap_{k\ge K}E_k,
 \qquad \PP(\Xi_{\rm typ})=1.
\end{equation}
Thus for every $\varpi\in\Xi_{\rm typ}$ there is a finite
$K(\varpi)$ such that, for every $k\ge K(\varpi)$,
\eqref{eq:random-forward} holds for the configuration at scale $a_k$.
Because $E_k$ depends only on the configuration and not on $P$, the same
realization and the same index $K(\varpi)$ work simultaneously for every
$P\ge P_0$ satisfying $Pd_k^{1/10}\le c_m$. In particular one may choose
a scale-dependent deterministic $P_k$ after the random realization has
been fixed, provided $P_kd_k^{1/10}\le c_m$.

Nothing special is used about the dyadic sequence except summability. If
$\sum_k d_k^m<\infty$, then
$\sum_k\PP(E_k^c)<\infty$ and the same Borel--Cantelli argument proves
\eqref{eq:random-typical-set}. No independence between different scales
is required. This is the precise almost-sure form of the random extension;
the analytic Foldy--Lax and boundary-map estimates remain deterministic
statements for each realization belonging to the corresponding event.
\section{Linearization of the N-D map and reconstruction}
\label{sec:linearization}

\subsection{The reference resolvent}
The reference operator for the prescribed density is
\[
 L_P=-\Delta+P^2\rho(x).
\]
Let $R_P$ be its inverse with zero Neumann condition. For $f\in\cX$,
denote by $v_P^f$ and $\Lambda_P^{\refe}$ the solution and the boundary
map of
\begin{equation}\label{eq:reference-PDE}
 L_Pv_P^f=0\quad\text{in }\Omega,\qquad
 \partial_\nu v_P^f=f,\qquad
 \Lambda_P^{\refe}f=\tr v_P^f.
\end{equation}
The operator $L_P$, including its variable coefficient $\rho$, is
known.

\begin{lemma}\label{lem:reference-resolvent}
For $P\ge P_0\ge1$,
\begin{align}
 \|R_P\|_{L^2\to L^2}&\le(\rho_*P^2)^{-1},&
 \|R_P\|_{L^2\to H^1}+\|\tr R_P\|_{L^2\to\cY}&\le CP^{-1},
 \label{eq:reference-resolvent}\\
 \|v_P^f\|_{H^1}&\le C\|f\|_{\cX},&
 \|v_P^f\|_{L^2}&\le CP^{-1}\|f\|_{\cX}.
 \label{eq:reference-data}
\end{align}
\end{lemma}
\begin{proof}
Let $z=R_PF$. Testing its equation with $\overline z$ gives
\[
 \|\nabla z\|_2^2+P^2\int_\Omega\rho|z|^2
 =\int_\Omega F\overline z.
\]
Dropping the gradient term and using Cauchy--Schwarz yields
$\|z\|_2\le(\rho_*P^2)^{-1}\|F\|_2$. Substitution into the same
identity gives $\|\nabla z\|_2\le CP^{-1}\|F\|_2$.
The trace theorem proves the remaining estimate in
\eqref{eq:reference-resolvent}.

For the reference solution, its weak equation and the trace theorem
give
\[
 c\|v_P^f\|_P^2
 \le\left|\pair{f}{\tr\overline{v_P^f}}\right|
 \le C\|f\|_{\cX}\|v_P^f\|_P.
\]
Hence $\|v_P^f\|_P\le C\|f\|_{\cX}$, which implies both estimates
in \eqref{eq:reference-data}.
\end{proof}

\subsection{The linearized Neumann-to-Dirichlet map}
Set $Q_\omega=\omega_{a,P}^2n^2$ and
\[
 \eta_P=\frac{q_*}{\rho_*P^2}.
\]
 Increase $P_0$ so that $\eta_P\le1/2$. The effective solution satisfies
\begin{equation}\label{eq:Born-LS}
 u_{\rho,P}^f=v_P^f+R_P(Q_\omega u_{\rho,P}^f),
 \qquad
 u_{\rho,P}^f=\sum_{j=0}^{\infty}(R_PQ_\omega)^jv_P^f
 \quad\text{in }L^2(\Omega).
\end{equation}
The series follows from
$\|R_PQ_\omega\|_{L^2\to L^2}\le\eta_P<1$.

\begin{theorem}\label{thm:linearization}
Define, for $b\in L^\infty(\Omega)$,
\begin{equation}\label{eq:linearized-map}
 \cL_P(b)f=\tr R_P(bv_P^f).
\end{equation}
Then
\begin{equation}\label{eq:linearization}
 \Lambda_{\rho,P}^{\eff}
 =\Lambda_P^{\refe}+\omega_{a,P}^2\cL_P(n^2)+\cR_{P,1},
 \qquad
 \|\cR_{P,1}\|_{\cL(\cX,\cY)}\le CP^{-4}.
\end{equation}
More generally, after retaining the terms of orders $1,\ldots,K$ in
\eqref{eq:Born-LS}, the boundary remainder satisfies
\begin{equation}\label{eq:Born-tail}
 \|\cR_{P,K}\|_{\cL(\cX,\cY)}
 \le CP^{-2}\frac{\eta_P^K}{1-\eta_P}
 \le C_KP^{-2K-2},\qquad K\ge1.
\end{equation}
The constants are uniform in the frequency and the geometric parameters
of Theorem~\ref{thm:main}.
\end{theorem}
\begin{proof}
For $j\ge1$, take the trace of the $j$th term in
\eqref{eq:Born-LS} by separating its first resolvent:
\[
 \tr(R_PQ_\omega)^jv_P^f
 =\tr R_PQ_\omega(R_PQ_\omega)^{j-1}v_P^f.
\]
Lemma~\ref{lem:reference-resolvent} gives
\begin{align*}
 \|\tr(R_PQ_\omega)^jv_P^f\|_{\cY}
 &\le\|\tr R_P\|_{L^2\to\cY}
       \|Q_\omega\|_\infty\eta_P^{j-1}\|v_P^f\|_2\\
 &\le CP^{-2}\eta_P^{j-1}\|f\|_{\cX}.
\end{align*}
Thus the series of traces converges in $\cY$, uniformly for
$\|f\|_{\cX}\le1$. Its sum agrees with the trace of the effective
solution: the sum of all nonzero-order terms is $R_PQ_\omega
u_{\rho,P}^f$, which belongs to $H^1(\Omega)$ by
\eqref{eq:reference-resolvent}. Summing the displayed estimate from
$j=K+1$ proves \eqref{eq:Born-tail}. The first term is
$\omega_{a,P}^2\cL_P(n^2)f$, which proves
\eqref{eq:linearization}.
\end{proof}

The map in \eqref{eq:linearized-map} is the Fr\'echet derivative, at
$b=0$, of the Neumann-to-Dirichlet map for $L_P-b$. This follows from
the same series, with $Q_\omega$ replaced by $b$ since the terms of degree at
least two are bounded by $CP^{-4}\|b\|_\infty^2$ on a fixed bounded
set of $b$, for sufficiently large $P$. Thus the linearization concerns
the unknown zeroth-order coefficient, while the large known potential
$P^2\rho$ is retained in the reference equation.

Suppose that a measured map satisfies
\begin{equation}\label{eq:measurement-noise}
 \|\Lambda_{a,P}^{\mathrm{meas}}-\Lambda_{a,P}\|_{\cL(\cX,\cY)}
 \le\delta_{a,P}.
\end{equation}
Theorems~\ref{thm:main} and \ref{thm:linearization} imply
\begin{equation}\label{eq:linear-data-error}
 \left\|\omega_{a,P}^{-2}
       (\Lambda_{a,P}^{\mathrm{meas}}-\Lambda_P^{\refe})
       -\cL_P(n^2)\right\|_{\cL(\cX,\cY)}
 \le C\bigl(E(a,P)+\delta_{a,P}+P^{-4}\bigr).
\end{equation}
The factor $\omega_{a,P}^{-2}$ is uniformly bounded because $I$ is a
compact subset of $(0,\infty)$.

For the random construction in Section~\ref{sec:random}, one can make
the homogenization error smaller than the linearization remainder.
Indeed, choose
\begin{equation}\label{eq:polynomial-regime}
 P=d^{-\sigma},\qquad 0<\sigma<\frac1{50}.
\end{equation}
Then $P d^{1/10}\to0$, and
\begin{align}
 P^4\bigl(P^2a^h+P^4d+P d^{1/10}\bigr)
 &=d^{3h/(1-h)-6\sigma}
      +d^{1-8\sigma}+d^{1/10-5\sigma}\longrightarrow0.
 \label{eq:linearization-compatible-rate}
\end{align}
All exponents are positive for $h>1/3$. Consequently the
physical-to-effective error is $o(P^{-4})$ on the common events. This
comparison concerns operator norms; the effect of CGO boundary inputs
is estimated separately below.

\subsection{CGO solutions for the known reference operator}
We recall the precise form of the CGO input which we use. Fix
$-1<\beta<0$ and write
\[
 L^2_\beta(\R^3)
 =\{u:(1+|x|^2)^{\beta/2}u\in L^2(\R^3)\}.
\]
For $\zeta\in\C^3$ with $\zeta\cdot\zeta=0$ and
$|\zeta|=\sqrt2s$, the Faddeev inverse $g_\zeta$ of
$-\Delta-2\zeta\cdot\nabla$ is a convolution operator and satisfies
\begin{equation}\label{eq:Faddeev-estimate}
 \|g_\zeta F\|_{L^2_\beta}\le\frac C s\|F\|_{L^2_{\beta+1}},
 \qquad
 \|\nabla g_\zeta F\|_{L^2_\beta}\le C\|F\|_{L^2_{\beta+1}},
 \quad s\ge1.
\end{equation}
These are the weighted resolvent estimates used in the construction of
\cite{SU}; see also \cite[Sections 2--3]{KR}. The local conjugated
estimates and the perturbative construction are discussed in
\cite[Section 3.2]{Salo}. We use the canonical convolution inverse in
\eqref{eq:Faddeev-estimate}, rather than an arbitrary local right
inverse, since the boundary construction in
Appendix~\ref{sec:correction} requires its fundamental-solution
property. Estimate \eqref{eq:Faddeev-estimate} is the standard CGO
resolvent theorem used as an input here. The dependence on $P$ is
obtained by the following argument.

Choose a known real extension $\widetilde\rho\in C_c^1(\R^3)$ of
$\rho$, with support in a fixed ball $B_0\supset\overline\Omega$,
and put $V_P=P^2\widetilde\rho$. Positivity of this extension outside
$\Omega$ is not needed.

\begin{lemma}\label{lem:reference-CGO}
If $s\ge C(1+P^2)$, there exists a solution
\begin{equation}\label{eq:reference-CGO}
 v_\zeta(x)=e^{\zeta\cdot x}(1+r_\zeta(x)),
 \qquad (-\Delta+V_P)v_\zeta=0\quad\text{in }\R^3,
\end{equation}
with $r_\zeta\in L^2_\beta(\R^3)$ and
\begin{equation}\label{eq:CGO-remainder}
 \|r_\zeta\|_{L^2_\beta}\le C\frac{P^2}s,
 \qquad \|\nabla r_\zeta\|_{L^2_\beta}\le CP^2.
\end{equation}
The solution is determined by the known potential $V_P$ and the
Faddeev normalization.
\end{lemma}
\begin{proof}
Multiplication by $V_P$, whose support is fixed, maps
$L^2_\beta$ into $L^2_{\beta+1}$ with norm at most $CP^2$.
Thus \eqref{eq:Faddeev-estimate} implies
\[
 \|g_\zeta V_P\|_{L^2_\beta\to L^2_\beta}
 \le C P^2/s.
\]
Choose the constant in $s\ge C(1+P^2)$ so that this norm is at most
$1/2$. Then
\begin{equation}\label{eq:CGO-integral}
 r_\zeta=-(I+g_\zeta V_P)^{-1}g_\zeta V_P
\end{equation}
is well-defined; in the last factor $V_P$ is viewed as a function.
Its norm in $L^2_{\beta+1}$ is at most $CP^2$, so the first estimate
in \eqref{eq:CGO-remainder} follows. Also
$r_\zeta=-g_\zeta[V_P(1+r_\zeta)]$ and
\[
 \|V_P(1+r_\zeta)\|_{L^2_{\beta+1}}
 \le CP^2(1+\|r_\zeta\|_{L^2_\beta})\le CP^2.
\]
The gradient estimate in \eqref{eq:Faddeev-estimate} gives the second
bound. Applying $-\Delta-2\zeta\cdot\nabla$ to
\eqref{eq:CGO-integral} proves the equation in
\eqref{eq:reference-CGO}. Interior elliptic regularity gives
$v_\zeta\in H^2_{\mathrm{loc}}(\R^3)$.
\end{proof}

For $\xi\in\R^3$, choose orthonormal vectors $e_1,e_2$ perpendicular
to $\xi$; when $\xi=0$, choose any two orthonormal vectors. If
$s>|\xi|/2$, set
\begin{align}
 \zeta_1&=s e_1+\ii\left(\frac\xi2+
                   \sqrt{s^2-|\xi|^2/4}\,e_2\right),\notag\\
 \zeta_2&=-s e_1+\ii\left(\frac\xi2-
                   \sqrt{s^2-|\xi|^2/4}\,e_2\right).
 \label{eq:CGO-phases}
\end{align}
A direct computation gives
\[
 \zeta_j\cdot\zeta_j=0,\qquad |\zeta_j|=\sqrt2s,
 \qquad \zeta_1+\zeta_2=\ii\xi.
\]
Let $v_j=v_{\zeta_j}$ and $f_j=\partial_\nu v_j|_{\partial\Omega}$.
These fluxes are known and belong to $\cX$. By uniqueness in
\eqref{eq:reference-PDE}, $v_j|_\Omega=v_P^{f_j}$.

\begin{proposition}\label{prop:linearized-reconstruction}
For $s\ge C(1+P^2+|\xi|)$,
\begin{equation}\label{eq:linearized-Fourier}
 \left|\pair{\cL_P(n^2)f_1}{f_2}
       -\int_\Omega n^2(x)e^{\ii\xi\cdot x}\dd x\right|
 \le C\left(\frac{P^2}s+\frac{P^4}{s^2}\right).
\end{equation}
In particular, the exact linearized map $\cL_P(n^2)$ determines
$n^2$ for every fixed $P\ge P_0$.
\end{proposition}
\begin{proof}
Put $W=R_P(n^2v_1)$. Its weak equation is
$L_PW=n^2v_1$ with zero Neumann data. Test it with $v_2$ and use the
reference equation for $v_2$. The symmetry of the bilinear form gives
\begin{equation}\label{eq:linearized-identity}
 \pair{\cL_P(n^2)f_1}{f_2}
 =\pair{\tr W}{\partial_\nu v_2}
 =\int_\Omega n^2v_1v_2\dd x.
\end{equation}
By \eqref{eq:CGO-phases},
\[
 v_1v_2=e^{\ii\xi\cdot x}
           (1+r_{\zeta_1}+r_{\zeta_2}+r_{\zeta_1}r_{\zeta_2}).
\]
The linear remainder terms have integrals bounded by
$C\|n^2\|_\infty|\Omega|^{1/2}P^2/s$.
Cauchy--Schwarz bounds the quadratic term by $CP^4/s^2$.
This proves \eqref{eq:linearized-Fourier}.

For fixed $P$ and $\xi$, let $s\to\infty$. We obtain the Fourier
transform, with positive exponential convention, of
$n^2\one_\Omega$ at every $\xi$. Fourier inversion in $L^2(\R^3)$
then determines $n^2$ in $\Omega$.
\end{proof}

\subsection{The boundary-data amplification}
Let $R_\Omega$ be such that $\overline\Omega\subset B(0,R_\Omega)$.
The CGO estimates and $|\operatorname{Re}\zeta_j|=s$ give
\begin{equation}\label{eq:CGO-H1}
 \|v_j\|_{H^1(\Omega)}
 \le C(1+s+P^2)e^{R_\Omega s}.
\end{equation}
Since $\Delta v_j=P^2\rho v_j$ in $\Omega$, the weak normal trace
estimate of Lemma~\ref{lem:normal-trace} yields
\begin{equation}\label{eq:CGO-flux-bound}
 \|f_j\|_{\cX}\le C(1+s+P^2)e^{R_\Omega s}.
\end{equation}
Here we used $\|v_j\|_2\le Ce^{R_\Omega s}$, which follows from
$P^2/s\le C$. Define
$\mathcal B_{P,s,\xi}=\|f_1\|_{\cX}\|f_2\|_{\cX}$.
For the approximation formed directly from the measured map,
\[
 F_{a,P}^{\mathrm{lin}}(\xi,s)
 =\omega_{a,P}^{-2}
   \pair{(\Lambda_{a,P}^{\mathrm{meas}}-\Lambda_P^{\refe})f_1}{f_2},
\]
we have the complete estimate
\begin{equation}\label{eq:physical-linearized-Fourier}
 \begin{split}
 \left|F_{a,P}^{\mathrm{lin}}(\xi,s)
        -\int_\Omega n^2e^{\ii\xi\cdot x}\dd x\right|
 \le {}&C\left(\frac{P^2}s+\frac{P^4}{s^2}\right)\\
 &+C\mathcal B_{P,s,\xi}
                 \bigl(E(a,P)+\delta_{a,P}+P^{-4}\bigr).
 \end{split}
\end{equation}
Indeed, apply \eqref{eq:linear-data-error} to $f_1,f_2$ and then use
\eqref{eq:linearized-Fourier}.

\begin{remark}\label{rem:linearized-versus-measured}
The operator estimate $\|\cR_{P,1}\|\le CP^{-4}$ does not necessarily
justify discarding the nonlinear remainder after inserting CGO data.
The factor $\mathcal B_{P,s,\xi}$ in
\eqref{eq:physical-linearized-Fourier} may grow exponentially in $s$.
Thus \eqref{eq:linearization-compatible-rate} makes the homogenization
error smaller than the linearization error in operator norm, but it is
not, on its own, a reconstruction theorem for the directly linearized
physical data. Proposition~\ref{prop:linearized-reconstruction} concerns
the exact linearized map. Appendix~\ref{sec:correction} provides a
boundary correction using the full map and proves a reconstruction
result with its scale and noise conditions.

The underlying issue is a coupling between linearization accuracy and
the amplification of boundary-map errors. Increasing $P$ reduces the
operator-level Born remainder, but the variable-density CGO construction
used here is available for $s\ge C(1+P^2+|\xi|)$. In particular,
\eqref{eq:CGO-flux-bound} gives the generic estimate
\begin{equation}\label{eq:Born-CGO-amplification}
 \left|\pair{\cR_{P,1}f_1}{f_2}\right|
 \le CP^{-4}\mathcal B_{P,s,\xi}
 \le CP^{-4}(1+s+P^2)^2e^{2R_\Omega s}.
\end{equation}
The last expression does not tend to zero as $P$ increases in the
admissible CGO regime. This shows a limitation of the operator-norm
estimate, namely it does not give a lower bound on the actual tested remainder.
The boundary correction removes the Born truncation, while
\eqref{eq:full-data-error} retains exponential amplification of the
physical approximation error and measurement noise. Thus $P$, $s$, and
the Fourier cutoff must be selected jointly relative to the available
data accuracy. Subsection~\ref{sec:Fourier-band} makes this selection
quantitative.

\end{remark}

When $\rho\equiv\rho_0$ is constant, the reference equation admits the
explicit solutions $e^{\zeta_j\cdot x}$ with
\[
 \zeta_1=\sqrt{P^2\rho_0+|\xi|^2/4}\,e_1+\ii\xi/2,
 \qquad
 \zeta_2=-\sqrt{P^2\rho_0+|\xi|^2/4}\,e_1+\ii\xi/2.
\]
Their product is exactly $e^{\ii\xi\cdot x}$, recovering the constant
reference construction of \cite[Remark 1.7]{JST}. For nonconstant
$\rho$, the remainders in Lemma~\ref{lem:reference-CGO} are required.
In both cases, the complex vectors are spatial phases of auxiliary
solutions and the incident frequency in \eqref{eq:tuning} remains real.
\appendix
\section{Reconstruction from the full boundary map}
\label{sec:correction}

This appendix gives a reconstruction procedure which uses the full
effective map, rather than replacing it by its first-order term.
It also transfers the procedure to the measured droplet map. The
boundary equation is a reference-potential version of the constructive
CGO equations in \cite{Nachman88,Novikov88}. We include its derivation
and the perturbation estimates, since their constants must be compared
with the error in Theorem~\ref{thm:main}.

\subsection{Dirichlet maps and reference fundamental solutions} \label{DM: Ref_FD}
We use Dirichlet-to-Neumann maps here because they enter the
boundary integral equation. In the coercive regime they are the inverses
of the Neumann-to-Dirichlet maps used above; see \eqref{eq:DN-bound}.
Fix the same real frequency $\omega=\omega_{a,P}$ and set
$Q=\omega^2n^2$. Denote the Dirichlet-to-Neumann maps for
$L_P-Q$ and $L_P$ on $\Omega$ by $\cD_P$ and $\cD_P^0$, respectively.
Both Dirichlet problems are coercive for $P\ge P_0$. Moreover,
\begin{equation}\label{eq:DN-bound}
 \|\cD_P\|_{\cL(\cY,\cX)}+
 \|\cD_P^0\|_{\cL(\cY,\cX)}\le C(1+P^2),
 \qquad
 \cD_P=(\Lambda_{\rho,P}^{\eff})^{-1},\quad
 \cD_P^0=(\Lambda_P^{\refe})^{-1}.
\end{equation}
To prove the norm bound, take a bounded trace extension
$E:\cY\to H^1(\Omega)$. The Dirichlet solution $u_h$ minimizes the
positive energy with trace $h$. Thus
\[
 \mathfrak a_\rho(u_h,u_h)
 \le\mathfrak a_\rho(Eh,Eh)
 \le C(1+P^2)\|h\|_{\cY}^2.
\]
For traces $h,k$, symmetry and the weak equations give
$\pair{\cD_Ph}{\overline k}=\mathfrak a_\rho(u_h,u_k)$.
Cauchy--Schwarz in this positive form proves the first bound in
\eqref{eq:DN-bound}; the proof for $\cD_P^0$ is the same. In particular,
\begin{equation}\label{eq:Dirichlet-H1}
 \|u_h\|_{H^1(\Omega)}+\|v_h\|_{H^1(\Omega)}
 \le C(1+P^2)^{1/2}\|h\|_{\cY},
\end{equation}
where $v_h$ is the reference Dirichlet solution. The inverse relations
in \eqref{eq:DN-bound} follow from uniqueness of both the Dirichlet
and Neumann problems. There is no restriction to mean-zero boundary
data, because the positive zeroth-order term removes the constant
Neumann kernel.

Retain the compactly supported extension $V_P=P^2\widetilde\rho$
from Section~\ref{sec:linearization}. On compactly supported volume
sources define
\begin{align}
 \cG_\zeta^0
 &=e^{\zeta\cdot x}g_\zeta e^{-\zeta\cdot x},\notag\\
 \cG_{P,\zeta}
 &=e^{\zeta\cdot x}(I+g_\zeta V_P)^{-1}
                              g_\zeta e^{-\zeta\cdot x}.
 \label{eq:reference-Faddeev}
\end{align}
Here $g_\zeta$ is the canonical inverse in
\eqref{eq:Faddeev-estimate}. The operator $\cG_\zeta^0$ is convolution
with a fundamental solution of $-\Delta$, whereas
$\cG_{P,\zeta}$ is a fundamental-solution operator for $-\Delta+V_P$.
Neither operator is the Neumann inverse $R_P$.
They satisfy the resolvent identity
\begin{equation}\label{eq:Faddeev-resolvent-identity}
 \cG_{P,\zeta}=\cG_\zeta^0
                         -\cG_{P,\zeta}V_P\cG_\zeta^0.
\end{equation}

For $\phi\in\cX$, let $\phi\delta_{\partial\Omega}$ denote the
surface distribution whose action on a test function is
$\pair{\phi}{\tr\varphi}$. Define the reference single-layer potential
and its boundary trace by
\begin{equation}\label{eq:reference-single-layer}
 \mathscr S_{P,\zeta}\phi
       =\cG_{P,\zeta}(\phi\delta_{\partial\Omega}),
 \qquad S_{P,\zeta}\phi=\tr\mathscr S_{P,\zeta}\phi.
\end{equation}
The meaning of \eqref{eq:reference-single-layer} for rough boundary
data, and the bounds required below, follow from the next lemma.

\begin{lemma}\label{lem:Faddeev-mapping}
Let $B_1$ be a fixed ball containing the support of $V_P$ and
$\overline\Omega$ in its interior. There are constants $C_1,C_2$,
independent of $P,s$ and of the direction of $\zeta$, such that, for
$s\ge C_1(1+P^2)$,
\begin{align}
 \|\cG_{P,\zeta}F\|_{H^2(B_1)}
       +\|\cG_\zeta^0F\|_{H^2(B_1)}
 &\le C_2e^{C_2s}\|F\|_{L^2(B_1)},
 \label{eq:volume-Faddeev-bound}\\
 \|\mathscr S_{P,\zeta}\phi\|_{H^1(B_1)}
       +\|S_{P,\zeta}\phi\|_{\cY}
 &\le C_2e^{C_2s}\|\phi\|_{\cX}.
 \label{eq:single-layer-bound}
\end{align}
In the first inequality, $F$ is extended by zero outside $B_1$.
For compactly supported $w\in H^1(\R^3)$, one has
\begin{equation}\label{eq:compact-left-inverse}
 \cG_{P,\zeta}(-\Delta+V_P)w=w.
\end{equation}
The identity is understood through the continuous extension of these
operators to the corresponding compactly supported distributions.
\end{lemma}
\begin{proof}
For a compactly supported source $F$, put
$z=(I+g_\zeta V_P)^{-1}g_\zeta F$.
The contraction estimate in Lemma~\ref{lem:reference-CGO} gives
$\|z\|_{L^2_\beta}\le C s^{-1}\|F\|_2$.
As $z=g_\zeta(F-V_Pz)$, \eqref{eq:Faddeev-estimate} also gives
$\|\nabla z\|_{L^2_\beta}\le C\|F\|_2$.
The equation
\[
 -\Delta z=2\zeta\cdot\nabla z+F-V_Pz
\]
and interior elliptic regularity on a slightly larger fixed ball imply
$\|z\|_{H^2(B_1)}\le C(1+s)\|F\|_2$.
Weights are equivalent to constants on each of these fixed balls.
Applying this estimate to $e^{-\zeta\cdot x}F$ and differentiating the
factor $e^{\zeta\cdot x}$ proves
\eqref{eq:volume-Faddeev-bound}, after absorbing polynomial factors
in $s$ into $e^{Cs}$. The free case is obtained by setting $V_P=0$.

For the boundary source, solve
\[
 (-\Delta+1)b=0\quad\text{in }\Omega,\qquad
 \partial_\nu b=\phi.
\]
Lax--Milgram gives $\|b\|_{H^1(\Omega)}\le C\|\phi\|_{\cX}$.
In distributions on $\R^3$,
\begin{equation}\label{eq:surface-source-decomposition}
 \phi\delta_{\partial\Omega}
 =-\operatorname{div}(\one_\Omega\nabla b)+\one_\Omega b.
\end{equation}
This follows by testing the weak Neumann equation with the restriction
of a smooth function on $\R^3$. Since the free operator is convolution,
it commutes with derivatives of compactly supported distributions.
Hence its potential of the source in
\eqref{eq:surface-source-decomposition} equals
\[
 \mathscr S_\zeta^0\phi
 =-\sum_{k=1}^3\partial_k\cG_\zeta^0
                         (\one_\Omega\partial_kb)
       +\cG_\zeta^0(\one_\Omega b).
\]
The free part of \eqref{eq:volume-Faddeev-bound} bounds this expression
in $H^1(B_1)$ by $Ce^{Cs}\|\phi\|_{\cX}$.
By \eqref{eq:Faddeev-resolvent-identity},
\[
 \mathscr S_{P,\zeta}\phi
 =\mathscr S_\zeta^0\phi
       -\cG_{P,\zeta}(V_P\mathscr S_\zeta^0\phi).
\]
The product on the right is a compactly supported $L^2$ source.
Another application of \eqref{eq:volume-Faddeev-bound}, using
$P^2\le Cs$, proves the $H^1$ estimate in
\eqref{eq:single-layer-bound}. Its boundary estimate follows from
the trace theorem. This argument defines the potential for every
$\phi\in\cX$ and agrees with the usual kernel definition on smooth
densities.

Finally, convolution with a fundamental solution gives
$\cG_\zeta^0(-\Delta)w=w$ for compactly supported smooth $w$.
Using \eqref{eq:Faddeev-resolvent-identity}, we obtain
\begin{align*}
 \cG_{P,\zeta}(-\Delta+V_P)w
 &=\cG_\zeta^0(-\Delta)w
       -\cG_{P,\zeta}V_P\cG_\zeta^0(-\Delta)w
       +\cG_{P,\zeta}V_Pw\\
 &=w.
\end{align*}
Approximation by compactly supported smooth functions, together with
the preceding volume and divergence-source bounds, extends this
identity to compactly supported $H^1$ functions. This proves
\eqref{eq:compact-left-inverse}.
\end{proof}

\subsection{The boundary integral equation}
Let $v_\zeta$ be the known reference CGO solution and let $u_\zeta$ be
the CGO solution for the full potential
$V_P-\one_\Omega Q$. The proof of
Lemma~\ref{lem:reference-CGO}, applied to this bounded compactly
supported potential, gives
\begin{equation}\label{eq:full-CGO}
 u_\zeta=e^{\zeta\cdot x}(1+r_\zeta^Q),\qquad
 \|r_\zeta^Q\|_{L^2_\beta}\le C\frac{1+P^2}s,
 \qquad s\ge C(1+P^2).
\end{equation}
Subtracting the conjugated reference and full equations gives
\begin{equation}\label{eq:full-CGO-volume}
 u_\zeta=v_\zeta+
              \cG_{P,\zeta}(\one_\Omega Q u_\zeta).
\end{equation}
Define
\begin{equation}\label{eq:boundary-B}
 \mathcal A_{P,\zeta}
 =I+S_{P,\zeta}(\cD_P-\cD_P^0):\cY\longrightarrow\cY.
\end{equation}

\begin{proposition}\label{prop:boundary-equation}
For $s\ge C(1+P^2)$, the operator
$\mathcal A_{P,\zeta}$ is invertible. The trace of the full CGO
solution is obtained from
\begin{equation}\label{eq:boundary-equation}
 \mathcal A_{P,\zeta}h_\zeta=\tr v_\zeta,
 \qquad h_\zeta=\tr u_\zeta.
\end{equation}
Furthermore,
\begin{equation}\label{eq:boundary-inverse-envelope}
 \|\mathcal A_{P,\zeta}^{-1}\|_{\cY\to\cY}\le Ce^{Cs}.
\end{equation}
\end{proposition}
\begin{proof}
We first derive the boundary identity for an arbitrary trace
$h\in\cY$. Let $u_h$ and $v_h$ be the full and reference Dirichlet
solutions with this trace. Set
\[
 w_h=u_h-v_h\in H_0^1(\Omega),\qquad
 \psi=(\cD_P-\cD_P^0)h.
\]
Extend $w_h$ by zero to $\R^3$, obtaining a compactly supported
$H^1$ function $\widetilde w_h$. The weak equations imply
\begin{equation}\label{eq:zero-extension-identity}
 (-\Delta+V_P)\widetilde w_h
   =\one_\Omega Q u_h+\psi\delta_{\partial\Omega}.
\end{equation}
In particular, the boundary term has a plus sign: testing the left side
and integrating the gradient term over $\Omega$ produces
$\pair{\partial_\nu w_h}{\tr\varphi}$.
Apply \eqref{eq:compact-left-inverse} to
\eqref{eq:zero-extension-identity}:
\begin{equation}\label{eq:zero-extension-potential}
 \widetilde w_h
 =\cG_{P,\zeta}(\one_\Omega Q u_h)
                  +\mathscr S_{P,\zeta}\psi.
\end{equation}
All terms have a common $H^1$ trace near $\partial\Omega$.
Since $\tr\widetilde w_h=0$, taking traces gives
\begin{equation}\label{eq:trace-volume-identity}
 \tr\cG_{P,\zeta}(\one_\Omega Q u_h)
      =-S_{P,\zeta}(\cD_P-\cD_P^0)h.
\end{equation}
The function
$u_h-\cG_{P,\zeta}(\one_\Omega Q u_h)$ solves the reference equation
in $\Omega$ and has trace $\mathcal A_{P,\zeta}h$.
For $u_h=u_\zeta|_\Omega$, \eqref{eq:full-CGO-volume} identifies this
function with $v_\zeta$, proving \eqref{eq:boundary-equation}.

We next prove bijectivity, without assuming a boundary Fredholm
alternative. Given $g\in\cY$, solve the reference Dirichlet problem
with trace $g$, and call its solution $v_g$. Seek $u$ on $\Omega$ from
\begin{equation}\label{eq:volume-for-general-boundary}
 u=v_g+\cG_{P,\zeta}(\one_\Omega Q u).
\end{equation}
After conjugation by $e^{-\zeta\cdot x}$, the operator on the right
has $L^2(\Omega)$ norm at most $C\|Q\|_\infty/s$.
This follows from \eqref{eq:Faddeev-estimate} and the inverse of
$I+g_\zeta V_P$; source and observation weights are bounded on $\Omega$.
For sufficiently large $s$, this norm is at most $1/2$.
Thus \eqref{eq:volume-for-general-boundary} has a unique solution.
Lemma~\ref{lem:Faddeev-mapping} then gives $u\in H^1(\Omega)$.
Applying the differential operator shows that $u$ solves the full
homogeneous equation. If $h=\tr u$, the preceding boundary identity
yields $\mathcal A_{P,\zeta}h=g$. This proves surjectivity.

If $\mathcal A_{P,\zeta}h=0$, the reference solution
$u_h-\cG_{P,\zeta}(\one_\Omega Q u_h)$ has zero Dirichlet trace,
and hence is zero. The same volume contraction implies $u_h=0$,
so $h=0$. This proves injectivity.

Finally, the volume contraction gives
$\|u\|_{L^2(\Omega)}\le Ce^{Cs}\|v_g\|_{L^2(\Omega)}$.
Use \eqref{eq:volume-for-general-boundary},
\eqref{eq:volume-Faddeev-bound} and
\eqref{eq:Dirichlet-H1}, and absorb powers of $P$ into the exponential
since $P^2\le Cs$. This yields
$\|u\|_{H^1(\Omega)}\le Ce^{Cs}\|g\|_{\cY}$.
Taking the trace proves \eqref{eq:boundary-inverse-envelope}.
\end{proof}

\subsection{The corrected Fourier identity}
Choose $\zeta_1,\zeta_2$ as in \eqref{eq:CGO-phases}. Let
$g_1=\tr v_{\zeta_1}$, $k_2=\tr v_{\zeta_2}$ and
$h_1=\mathcal A_{P,\zeta_1}^{-1}g_1$. Define
\begin{equation}\label{eq:corrected-Fourier}
 F_P(\xi,s)=\omega^{-2}
                   \pair{(\cD_P^0-\cD_P)h_1}{k_2}.
\end{equation}
Green's identity gives
\begin{equation}\label{eq:corrected-integral}
 F_P(\xi,s)=\int_\Omega n^2 u_{\zeta_1}v_{\zeta_2}\dd x.
\end{equation}
Indeed the full weak equation gives the reference bilinear energy of
$u_{\zeta_1},v_{\zeta_2}$ minus $\int Q u_{\zeta_1}v_{\zeta_2}$
as $\pair{\cD_Ph_1}{k_2}$, while the reference equation gives that
energy as $\pair{\cD_P^0h_1}{k_2}$. Subtraction proves
\eqref{eq:corrected-integral}.
Using \eqref{eq:full-CGO} and \eqref{eq:CGO-remainder} in the product,
as in Proposition~\ref{prop:linearized-reconstruction}, gives
\begin{equation}\label{eq:corrected-Fourier-error}
 \left|F_P(\xi,s)-\int_\Omega n^2e^{\ii\xi\cdot x}\dd x\right|
 \le C\frac{1+P^2}s,
 \qquad s\ge C(1+P^2+|\xi|).
\end{equation}
There is no Born remainder in \eqref{eq:corrected-Fourier-error} and the
unknown full CGO trace is obtained through the boundary equation.

\subsection{Perturbation of the measured map}
Put
\[
 \eta_N=\|\Lambda_{a,P}^{\mathrm{meas}}-
                         \Lambda_{\rho,P}^{\eff}\|_{\cL(\cX,\cY)}.
\]
On $\cE_a$, Theorem~\ref{thm:main} and
\eqref{eq:measurement-noise} give
$\eta_N\le CE(a,P)+\delta_{a,P}$. If
$\|\cD_P\|\eta_N\le1/2$, then the measured map is invertible and
\begin{equation}\label{eq:measured-DN}
 \cD_P^{\mathrm{meas}}=(\Lambda_{a,P}^{\mathrm{meas}})^{-1},\qquad
 \eta_D:=\|\cD_P^{\mathrm{meas}}-\cD_P\|
 \le2\|\cD_P\|^2\eta_N
 \le C(1+P^2)^2\eta_N.
\end{equation}
For example, this follows by writing
$\Lambda_{a,P}^{\mathrm{meas}}
=\Lambda_{\rho,P}^{\eff}[I+\cD_P
(\Lambda_{a,P}^{\mathrm{meas}}-\Lambda_{\rho,P}^{\eff})]$ and
inverting the second factor by a Neumann series.

For the first phase, abbreviate
\[
 \kappa=\|\mathcal A_{P,\zeta_1}^{-1}\|,\qquad
 \sigma_S=\|S_{P,\zeta_1}\|,\qquad
 \Delta_D=\|\cD_P-\cD_P^0\|.
\]
Replace $\cD_P$ by $\cD_P^{\mathrm{meas}}$ in
\eqref{eq:boundary-B}, obtaining $\mathcal A^{\mathrm{meas}}$.
If $\kappa\sigma_S\eta_D\le1/2$, it is invertible and, for
$h_1^{\mathrm{meas}}=(\mathcal A^{\mathrm{meas}})^{-1}g_1$,
\begin{align}
 \|h_1^{\mathrm{meas}}\|_{\cY}&\le2\kappa\|g_1\|_{\cY},\notag\\
 \|h_1^{\mathrm{meas}}-h_1\|_{\cY}
 &\le2\kappa^2\sigma_S\eta_D\|g_1\|_{\cY}.
 \label{eq:trace-perturbation}
\end{align}
The second estimate follows by subtracting the two boundary equations:
$\mathcal A(h_1^{\mathrm{meas}}-h_1)
=-S_{P,\zeta_1}(\cD_P^{\mathrm{meas}}-\cD_P)h_1^{\mathrm{meas}}$.

Define the measured version of \eqref{eq:corrected-Fourier} by
\[
 F_{a,P}^{\mathrm{meas}}(\xi,s)
 =\omega^{-2}\pair{(\cD_P^0-\cD_P^{\mathrm{meas}})
                              h_1^{\mathrm{meas}}}{k_2}.
\]
Subtracting \eqref{eq:corrected-Fourier} and using
\eqref{eq:trace-perturbation} yields the explicit data estimate
\begin{equation}\label{eq:conditioning-explicit}
 \begin{split}
 \left|F_{a,P}^{\mathrm{meas}}(\xi,s)
       -\int_\Omega n^2e^{\ii\xi\cdot x}\dd x\right|
 \le {}&C\frac{1+P^2}s\\
 &+2\omega^{-2}\kappa\|g_1\|_{\cY}\|k_2\|_{\cY}
                  (1+\Delta_D\kappa\sigma_S)\eta_D.
 \end{split}
\end{equation}
No conditioning factor has been included in an unspecified
$P$-dependent constant in this inequality.

\begin{theorem}\label{thm:physical-CGO}
There are constants $C,c>0$ such that, on $\cE_a$, the measured
boundary correction is well-defined if \eqref{eq:joint-regime} holds,
$s\ge C(1+P^2+|\xi|)$, and
\begin{equation}\label{eq:CGO-data-smallness}
 (1+P^2)^2e^{Cs}\bigl(E(a,P)+\delta_{a,P}\bigr)\le c.
\end{equation}
In that case,
\begin{equation}\label{eq:full-data-error}
 \left|F_{a,P}^{\mathrm{meas}}(\xi,s)
       -\int_\Omega n^2e^{\ii\xi\cdot x}\dd x\right|
 \le C\frac{1+P^2}s
    +C(1+P^2)^2e^{Cs}\bigl(E(a,P)+\delta_{a,P}\bigr).
\end{equation}
The constants are independent of $a,P,s,\xi$ and the realization in
the stated ranges.
\end{theorem}
\begin{proof}
Lemmas~\ref{lem:reference-CGO} and \ref{lem:Faddeev-mapping}, the
trace theorem, and Proposition~\ref{prop:boundary-equation} give
\[
 \kappa+\sigma_S+\|g_1\|_{\cY}+\|k_2\|_{\cY}\le Ce^{Cs}.
\]
Also $\Delta_D\le C(1+P^2)$ by \eqref{eq:DN-bound}.
Since $P^2\le Cs$, products of these quantities are bounded by
$Ce^{C's}$ after increasing the constant in the exponent.
Condition \eqref{eq:CGO-data-smallness}, with a sufficiently large
$C$ and sufficiently small $c$, implies both inverse conditions
preceding \eqref{eq:measured-DN} and
\eqref{eq:trace-perturbation}. Substitution of
\eqref{eq:measured-DN} into \eqref{eq:conditioning-explicit} proves
\eqref{eq:full-data-error}. All estimates for the phases are uniform
in their directions, so the same constants apply to all $\xi$ in the
stated range.
\end{proof}

\subsection{Quantitative trade-off and recoverable Fourier band}
\label{sec:Fourier-band}
We extract a finite-accuracy criterion from
Theorem~\ref{thm:physical-CGO}. Put
$M_P=1+P^2$ and $\eta_{a,P}=E(a,P)+\delta_{a,P}$, and, for $K\ge0$, define
\[
 \mathcal E_F(a,P,s,K)
 =\sup_{|\xi|\le K}
 \left|F_{a,P}^{\mathrm{meas}}(\xi,s)
       -\int_\Omega n^2(x)e^{\ii\xi\cdot x}\dd x\right|.
\]
Since the estimates are uniform in the directions of the CGO phases,
there are constants $C_0,C_1,C_2,C_{\mathrm e},c_{\mathrm{inv}}>0$,
independent of the parameters and the realization, such that, on
$\cE_a$ and under \eqref{eq:joint-regime}, the boundary correction is
well-defined throughout $|\xi|\le K$ whenever
\begin{equation}\label{eq:band-admissibility}
 s\ge C_0(M_P+K),\qquad
 M_P^2 e^{C_{\mathrm e}s}\eta_{a,P}\le c_{\mathrm{inv}}.
\end{equation}
In this range,
\begin{equation}\label{eq:band-error-envelope}
 \mathcal E_F(a,P,s,K)
 \le C_1\frac{M_P}{s}
      +C_2M_P^2e^{C_{\mathrm e}s}\eta_{a,P}.
\end{equation}

\begin{corollary}[A sufficient Fourier-band criterion]
\label{cor:Fourier-band}
Assume \eqref{eq:joint-regime} and $\eta_{a,P}>0$. For a prescribed
Fourier error tolerance $\tau>0$, set
\begin{equation}\label{eq:band-noise-budget}
 m_\tau=\min\left\{c_{\mathrm{inv}},\frac{\tau}{2C_2}\right\},
 \qquad
 S_\tau(a,P)=\frac{1}{C_{\mathrm e}}
       \log\frac{m_\tau}{M_P^2\eta_{a,P}}.
\end{equation}
If the interval
\begin{equation}\label{eq:CGO-window}
 \max\left\{C_0(M_P+K),\frac{2C_1M_P}{\tau}\right\}
 \le s\le S_\tau(a,P)
\end{equation}
is nonempty, then every $s$ in this interval gives a well-defined
measured boundary correction on $\cE_a$ and
$\mathcal E_F(a,P,s,K)\le\tau$.
\end{corollary}
\begin{proof}
The lower bound in \eqref{eq:CGO-window} gives CGO admissibility and
$C_1M_P/s\le\tau/2$. The upper bound gives
$M_P^2e^{C_{\mathrm e}s}\eta_{a,P}\le m_\tau$, which ensures the
perturbation condition in \eqref{eq:band-admissibility} and bounds
the second term of \eqref{eq:band-error-envelope} by $\tau/2$.
\end{proof}

Consequently, if
\begin{equation}\label{eq:band-feasibility}
 S_\tau(a,P)\ge
 \max\left\{C_0M_P,\frac{2C_1M_P}{\tau}\right\},
\end{equation}
every Fourier radius in the range
\begin{equation}\label{eq:certified-Fourier-radius}
 0\le K\le K_{\mathrm{cert}}(a,P,\tau)
 :=\frac{1}{C_0C_{\mathrm e}}
       \log\frac{m_\tau}{(1+P^2)^2\eta_{a,P}}-(1+P^2)
\end{equation}
is guaranteed at tolerance $\tau$. The feasibility condition is part
of this conclusion. If $\eta_{a,P}=0$, there is no upper restriction
on $s$ in this criterion. These are sufficient bounds for
certified Fourier accuracy.

The window \eqref{eq:CGO-window} displays the parameter trade-off:
the reference potential, Fourier radius, and desired accuracy impose
its lower endpoint, while the physical and measurement errors impose
its upper endpoint. Equivalently, \eqref{eq:CGO-window}
is nonempty exactly when
\[
 \eta_{a,P}\le \frac{m_\tau}{M_P^2}
 \exp\!\left[-C_{\mathrm e}
 \max\!\left\{C_0(M_P+K),\frac{2C_1M_P}{\tau}\right\}\right].
\]
This explicit compatibility condition is achieved in the joint regime
of Subsection~\ref{sec:joint-reconstruction}.

For a joint limit, write
$L_{a,P}=\log\bigl(1/(M_P^2\eta_{a,P})\bigr)\to\infty$.
Choosing $s=\theta L_{a,P}/C_{\mathrm e}$, with fixed $0<\theta<1$,
gives
\begin{equation}\label{eq:logarithmic-band-balance}
 \begin{split}
 \mathcal E_F(a,P,s,K)
 &\le \frac{C_1C_{\mathrm e}}{\theta}
             \frac{M_P}{L_{a,P}}
       +C_2e^{-(1-\theta)L_{a,P}},\\
 K&\le\frac{\theta L_{a,P}}{C_0C_{\mathrm e}}-M_P,
 \end{split}
\end{equation}
for all sufficiently large $L_{a,P}$ in the physical admissible regime.
Thus $M_P=o(L_{a,P})$ permits vanishing uniform Fourier error on an
expanding band. This is a sufficient compatibility condition for the
present estimates. Convergence of the reconstructed coefficient also
requires control of the accumulated error over that band, as shown
in Subsection \ref{sec:joint-reconstruction}.

\subsection{A sufficient joint reconstruction regime} \label{sec:joint-reconstruction} 
We give one choice of parameters which meets all the preceding
conditions. For the random construction of Section~\ref{sec:random},
fix $\nu>0$ and take, along integers $k\to\infty$,
\begin{equation}\label{eq:superfine-regime}
 P_k=k+P_0,\qquad s_k=P_k^{2+\nu},\qquad
 a_k=\exp(-P_k^{2+2\nu}),\qquad
 \delta_{a_k,P_k}\le\exp(-P_k^{2+2\nu}).
\end{equation}
Then $d_k=a_k^{(1-h)/3}$ and $P_kd_k^{1/10}\to0$.
On the events of Corollary~\ref{cor:random-forward}, the sum of the
forward and measurement errors is bounded by a polynomial in $P_k$
times $\exp(-c_hP_k^{2+2\nu})$, where one may take
\[
 c_h=\min\{h,(1-h)/3,(1-h)/30\}>0.
\]
This decay dominates the factor $e^{Cs_k}$. Thus
\eqref{eq:CGO-data-smallness} holds for large $k$, and
\begin{equation}\label{eq:uniform-frequency-recovery}
 \sup_{|\xi|\le R_k}
 \left|F_{a_k,P_k}^{\mathrm{meas}}(\xi,s_k)
       -\int_\Omega n^2e^{\ii\xi\cdot x}\dd x\right|
 \le CP_k^{-\nu},\qquad R_k=P_k^{\nu/3}.
\end{equation}
One can choose the orthonormal vectors in \eqref{eq:CGO-phases}
measurably in $\xi$, using finitely many coordinate choices. This makes
the reconstructed Fourier data measurable.

Define
\[
 n_k^2(x)=(2\pi)^{-3}\int_{|\xi|\le R_k}
                e^{-\ii\xi\cdot x}
                F_{a_k,P_k}^{\mathrm{meas}}(\xi,s_k)\dd\xi.
\]
By Plancherel's theorem and \eqref{eq:uniform-frequency-recovery},
\begin{align*}
 \|n_k^2-n^2\one_\Omega\|_{L^2(\R^3)}
 &\le C R_k^{3/2}P_k^{-\nu}
   +C\left\|\one_{\{|\xi|>R_k\}}
                   \widehat{n^2\one_\Omega}\right\|_{L^2(\R^3)}\\
 &\longrightarrow0.
\end{align*}
The first term is $CP_k^{-\nu/2}$ and the second is a Fourier tail of
an $L^2$ function. The exceptional probabilities $d_k^m$ are summable.
By Corollary~\ref{cor:probability}, this reconstruction holds for every
realization in a set of probability one, for all sufficiently large $k$.
The physical frequencies throughout this procedure are still those in
\eqref{eq:tuning}; only the auxiliary CGO phases and the spatial
resolution parameters tend to infinity.

\section{Auxiliary estimates}
\label{sec:auxiliary}

\subsection{Uniform background bounds}
We give the details which connect the background spectral assumption to
the kernel and solution estimates used in the point-interaction proof.
The underlying elliptic regularity and trace results can be found in
\cite{GT,McLean}; for Green kernels and boundary potential theory, see
also \cite{MitreaTaylor}.

\begin{lemma}\label{lem:background-bounds}
Under \eqref{eq:background-gap},
\begin{equation}\label{eq:background-H1}
 \|S^f\|_{H^1(\Omega)}\le C\|f\|_{\cX},
 \qquad \omega\in I,
\end{equation}
and \eqref{eq:green-bounds}--\eqref{eq:H-bounds} hold with uniform
constants.
\end{lemma}
\begin{proof}
Let $A_+=-\Delta+1$ with zero Neumann condition, let $R_+=A_+^{-1}$,
and put $W_\omega=1+q_\omega$. The operator $W_\omega$ is multiplication
by a uniformly bounded Lipschitz function, and
$A_{0,\omega}=A_+-W_\omega$.

First solve $A_+S_+^f=0$ with flux $f$. By Lax--Milgram and the trace
theorem, $\|S_+^f\|_{H^1}\le C\|f\|_{\cX}$. The background solution
has the representation
\[
 S^f=S_+^f+A_{0,\omega}^{-1}(W_\omega S_+^f).
\]
The second term has zero Neumann data and cancels
$A_{0,\omega}S_+^f=-W_\omega S_+^f$.
Global Neumann regularity gives
\[
 \|A_{0,\omega}^{-1}F\|_{H^2}
 \le C\bigl(\|F\|_2+\|A_{0,\omega}^{-1}F\|_2\bigr)
 \le C\|F\|_2,
\]
where the final constant is uniform by
\eqref{eq:background-gap}. This proves
\eqref{eq:background-H1}.

We next prove the Green-kernel value estimate without assuming that
$A_{0,\omega}$ is positive. The positive Neumann kernel $G_+$ of
$A_+$ has the standard bound
\[
 |G_+(x,y)|\le C|x-y|^{-1},\qquad
 \sup_{x\in\Omega}\|G_+(x,\cdot)\|_{L^2(\Omega)}\le C.
\]
The second bound follows from the first by integration in dimension
three. The two resolvent identities imply
\begin{equation}\label{eq:twice-resolvent}
 A_{0,\omega}^{-1}
 =R_++R_+W_\omega R_+
      +R_+W_\omega A_{0,\omega}^{-1}W_\omega R_+.
\end{equation}
For fixed $x,y$, the kernel of the second term has modulus at most
\[
 \|W_\omega\|_\infty
       \|G_+(x,\cdot)\|_2\|G_+(\cdot,y)\|_2\le C.
\]
The kernel of the third term is the $L^2$ pairing of
$W_\omega G_+(x,\cdot)$ with
$A_{0,\omega}^{-1}(W_\omega G_+(\cdot,y))$ and is also bounded by
$C$, using \eqref{eq:background-gap}. Thus
$|G_\omega(x,y)-G_+(x,y)|\le C$.
Since $\Omega$ is bounded, this proves
$|G_\omega(x,y)|\le C|x-y|^{-1}$.
The kernel is real and symmetric because the Neumann resolvent is
real and self-adjoint.

To obtain the gradient estimate, take $r=|x-y|$ and apply a scaled
interior gradient estimate, or a homogeneous Neumann boundary gradient
estimate, on a neighborhood of $x$ of radius comparable to $r$ which
does not contain $y$. On that neighborhood the kernel solves the
homogeneous equation, and its supremum is at most $C/r$.
Consequently $|\nabla_xG_\omega(x,y)|\le C/r^2$.
For radii larger than the size of a boundary chart, use a fixed chart
radius and then adjust the constant, since $r\le\operatorname{diam}\Omega$.
The estimates are uniform for the bounded family $q_\omega$ and for
the finite collection of $C^{2,\alpha}$ boundary charts. Symmetry gives
the same bound for $\nabla_yG_\omega$. This proves
\eqref{eq:green-bounds}.

It remains to estimate the regular part near the diagonal. Fix $y$
at distance at least $cb$ from $\partial\Omega$, and choose a ball
$B(y,R)\Subset\Omega$ with $R$ a fixed sufficiently small multiple of
$b$. The local parametrix
\[
 K_y(x)=\frac1{4\pi|x-y|}
                    -\frac{q_\omega(y)}{8\pi}|x-y|
\]
satisfies, distributionally,
\begin{align*}
 A_{0,\omega}K_y
 =\delta_y+[q_\omega(y)-q_\omega(x)]\Phi_0(x,y)
      +\frac{q_\omega(x)q_\omega(y)}{8\pi}|x-y|.
\end{align*}
Here we used $\Delta|x-y|=2/|x-y|$ away from $y$; this identity also
holds distributionally without a point mass. Therefore
$F_y=G_\omega(\cdot,y)-K_y$ solves
\begin{equation}\label{eq:parametrix-remainder}
 A_{0,\omega}F_y
 =[q_\omega(x)-q_\omega(y)]\Phi_0(x,y)
      -\frac{q_\omega(x)q_\omega(y)}{8\pi}|x-y|.
\end{equation}
The right-hand side is uniformly bounded because $q_\omega$ is
Lipschitz. Initially $F_y\in L^2$ locally, and interior elliptic
regularity applied to this distributional equation gives
$F_y\in H^2_{\mathrm{loc}}$, so that it is continuous near $y$.
On $\partial B(y,R)$, \eqref{eq:green-bounds} gives
$|F_y|\le C/R$.

Rescale the Dirichlet problem on this ball. For $R$ sufficiently small,
the norm of the zero-Dirichlet Newtonian solution operator followed by
multiplication by $q_\omega$ is at most $C R^2\|q_\omega\|_\infty<1/2$
on $L^\infty$. Splitting into the harmonic extension of the boundary
values and a zero-Dirichlet correction in
\eqref{eq:parametrix-remainder} gives
\[
 \|F_y\|_{L^\infty(B(y,R))}\le C(R^{-1}+R^2)\le C/R.
\]
The scaled interior gradient estimate then yields
$\|\nabla F_y\|_{L^\infty(B(y,R/2))}\le C/R^2$.
Since
\[
 H_\omega(x,y)=F_y(x)-\frac{q_\omega(y)}{8\pi}|x-y|,
\]
we obtain the asserted value and Lipschitz bounds when
$|x-y|<R/2$. For $|x-y|\ge R/2$ in the region specified in
\eqref{eq:H-bounds}, the estimates follow directly by subtracting
$\Phi_0$ in \eqref{eq:green-bounds}. Finally, symmetry gives the same
Lipschitz estimate in the second variable. This proves
\eqref{eq:H-bounds}, including the existence of the diagonal value.
\end{proof}

\subsection{A block-operator estimate}
\begin{lemma}\label{lem:block-schur}
Let $\mathcal V$ be a Hilbert space and let $A_{ij}:\mathcal V\to
\mathcal V$ satisfy $\|A_{ij}\|\le c_{ij}$ with $c_{ij}\ge0$.
Suppose
\[
 \sup_i\sum_jc_{ij}\le R,\qquad
 \sup_j\sum_ic_{ij}\le C_0.
\]
On the product space with norm
$\|v\|_w^2=w\sum_j\|v_j\|_{\mathcal V}^2$, one has
\[
 \left\|\left(\sum_jA_{ij}v_j\right)_i\right\|_w
 \le\sqrt{RC_0}\,\|v\|_w.
\]
\end{lemma}
\begin{proof}
For each row, the triangle inequality and weighted Cauchy--Schwarz give
\[
 \left\|\sum_jA_{ij}v_j\right\|_{\mathcal V}^2
 \le\left(\sum_jc_{ij}\|v_j\|_{\mathcal V}\right)^2
 \le\left(\sum_jc_{ij}\right)
                    \sum_jc_{ij}\|v_j\|_{\mathcal V}^2.
\]
Multiply by $w$, sum in $i$, and interchange the finite sums.
The result is at most
$RC_0w\sum_j\|v_j\|_{\mathcal V}^2$.
Taking square roots proves the assertion.
\end{proof}

\subsection{The weak normal trace}
\begin{lemma}\label{lem:normal-trace}
If $u\in H^1(\Omega)$ and $\Delta u\in L^2(\Omega)$, its weak normal
trace belongs to $\cX$ and satisfies
\begin{equation}\label{eq:weak-normal-trace}
 \|\partial_\nu u\|_{\cX}
 \le C\bigl(\|u\|_{H^1(\Omega)}+\|\Delta u\|_{L^2(\Omega)}\bigr).
\end{equation}
\end{lemma}
\begin{proof}
For $h\in\cY$, choose a bounded $H^1$ extension $Eh$ of the trace
and define
\[
 \pair{\partial_\nu u}{h}
 =\int_\Omega\nabla u\cdot\nabla(Eh)+\Delta u\,Eh.
\]
The definition is independent of the extension: the difference of two
extensions belongs to $H_0^1(\Omega)$, and the distributional definition
of $\Delta u$ makes the corresponding integral vanish.
Cauchy--Schwarz and $\|Eh\|_{H^1}\le C\|h\|_{\cY}$ prove
\eqref{eq:weak-normal-trace}. For smooth functions it agrees with the
classical normal derivative, by integration by parts.
\end{proof}

\end{document}